\input amstex
\documentstyle {amsppt}
\UseAMSsymbols \vsize 18cm \widestnumber\key{ZZZZZ}

\catcode`\@=11
\def\displaylinesno #1{\displ@y\halign{
\hbox to\displaywidth{$\@lign\hfil\displaystyle##\hfil$}&
\llap{$##$}\crcr#1\crcr}}
\def\ldisplaylinesno #1{\displ@y\halign{
\hbox to\displaywidth{$\@lign\hfil\displaystyle##\hfil$}&
\kern-\displaywidth\rlap{$##$} \tabskip\displaywidth\crcr#1\crcr}}
\catcode`\@=12

\refstyle{A}

\let \ol=\overline
\let \ul=\underline

\font\sit=cmti10 at 7.5pt

\font\nr=eufb7 at 10pt

 \font\srm=cmr10 at 7.5pt

\font\main=cmsy10 at 10pt

\font\smain=cmsy10 at 7.5pt \font\ssmain=cmsy10 at 5.625pt

\font \fin=lasy8 at 15.4 pt
\def \X{\mathop{\hbox{\nr X}^{\hbox{\srm nr}}}\nolimits}

\def \o{\mathop{\hbox{\main O}}\nolimits}
\def \so{\mathop{\hbox{\smain O}}\nolimits}
\def \sso{\mathop{\hbox{\ssmain O}}\nolimits}

\def \C{\mathop{\hbox{\main C}}\nolimits}

\def \A{\mathop{\hbox{\main A}}\nolimits}
\def \H{\mathop{\hbox{\main H}}\nolimits}

\def \B{\mathop{\hbox{\main B}}\nolimits}

\def \J{\mathop{\hbox{\main J}}\nolimits}
\def \R{\mathop{\hbox{\main R}}\nolimits}
\def \S{\mathop{\hbox{\main S}}\nolimits}

\def \sJ{\mathop{\hbox{\smain J}}\nolimits}

\def \sS{\mathop{\hbox{\smain S}}\nolimits}

\def \End{\mathop{\hbox{\rm End}}\nolimits}

\def \Hom{\mathop{\hbox{\rm Hom}}\nolimits}

\def \GL{\mathop{\hbox{\rm GL}}\nolimits}

\def \det{\mathop{\hbox{\rm det}}\nolimits}
\topmatter
\title Local Langlands correspondence for classical groups and affine Hecke algebras\endtitle

\rightheadtext{Classical groups and affine Hecke algebras}
\author Volker Heiermann  \endauthor
\address Aix Marseille Universit\'e, CNRS, Centrale Marseille, I2M UMR 7373, 13453, Marseille, France
\endaddress
\email volker.heiermann\@univ-amu.fr \endemail

\thanks The author has benefitted from a grant of Agence
Nationale de la Recherche with reference ANR-08-BLAN-0259-02 and with reference ANR-13-BS01-0012 FERPLAY.
\endthanks

\abstract Using the results of J. Arthur on the representation theory of classical groups with additional work by Colette Moeglin and its relation with representations
of affine Hecke algebras established by the author, we show that the category of smooth complex representations of a quasi-split $p$-adic classical group and its pure inner
forms is naturally decomposed into subcategories which are equivalent to a tensor product of categories of unipotent representations of classical groups. A statement of this kind had been conjectured by G. Lusztig. All classical groups (general linear, orthogonal, symplectic and unitary groups) appear in this context. We get also parameterizations of representations of affine Hecke algebras, which seem not all to be in the literature yet. All this sheds some light on what is known as the stable Bernstein center.

\endabstract

\endtopmatter
\document
Let $F$ be a non-Archimedean local field of characteristic $0$ - this assumption on the characteristic is used in \cite{A, M} but not in \cite{H3} -, and $n\geq 1$ an integer. The symbol $G$ will denote the group of $F$-rational points of a quasi-split classical group $\underline{G}$ of semi-simple rank $n$ defined over $F$. We will mean by that either a symplectic group or a (at least in the even rank case, the full, i.e. non connected) orthogonal group. (The case of unitary groups will be treated in the annex.) If $G$ is orthogonal, we will denote by $G^-$ its unique pure inner form [V, GGP]. If $G$ is symplectic, we will leave $G^-$ undefined (there is no pure inner form $\ne G$). We will often write $G^+$ for $G$ and denote by $Rep(G^\pm)$ the category of smooth complex representations of $G^\pm$.

J. Arthur has determined in \cite{A} the parameters of discrete and tempered $L$-packets of $G$, including the description of the $R$-groups. (The parametrization for the inner forms of $G$ is forthcoming.) C. Moeglin has deduced from this in \cite{M} the Langlands-Deligne parameters which correspond to supercuspidal representations (for both $G$ and $G^-$), including information on reducibility points. The author has used this information in \cite{H4} to deduce the parameters of the affine Hecke algebras that have been shown in \cite{H3} to correspond to the Bernstein components of the category $Rep(G^\pm)$.

The aim of the present work is to show that, putting together different Bernstein components, one obtains a natural decomposition of $Rep(G^\pm)$ into subcategories $\R^{\varphi _0}_F(\underline{G})$ which are equivalent to a tensor product of categories of unipotent representations of classical groups (in the sense of \cite{L2}).  (Because of the categorical nature of the tensor product, we have in fact to restrict to the full subcategory of finitely generated representations.) A statement of this kind had been conjectured by G. Lusztig \cite{L4, section 19}. The $\varphi _0$ can be seen as inertial classes of Langlands parameters for $G$ (i.e. modulo restriction to the inertial subgroup). When $^LG$ is orthogonal, a quasi-split outer form (equal to $G$ if $G$ is symplectic) will be added to obtain uniform statements. All classical groups (general linear, orthogonal, symplectic and unitary groups) appear in this context. Once the results in \cite{L2} appropriately generalized to symplectic, unitary and the (full) even orthogonal group, one should be able to compute multiplicities in standard modules from intersection cohomology as described in \cite{L4, section 19}.

Taking into account the local Langlands correspondence (which is in good shape now after the above mentioned results of Arthur and also of Moeglin - however no final account has been written yet, additional results on the preservation of local constants related to non standard $L$-functions are in \cite{He, CST}), we get from this  parameterizations of representations of affine Hecke algebras, which seem not all to be in the literature yet. In addition, we explain, how to get a direct correspondence for the irreducible representations in $\R^{\varphi _0}_F(\underline{G})$ by conjugacy classes of parameters $(s,u,\Xi )$ in a given complex reductive group (where $s$ is a semi-simple element, $u$ a unipotent element such that $sus^{-1}=u^q$ and $\Xi $ an irreducible representation of the group of components of the common centralizer of $s$ and $u$).

The plan of this paper is the following: in section {\bf 1.}, we summarize the results of Moeglin (based on Arthur's work) on the Langlands correspondence for supercuspidal representations of $G$. We recall the author's results relating the Bernstein components of $Rep(G^\pm)$ to representations of affine Hecke algebras and  give the definition of the categories $\R^{\varphi _0}_F(\underline{G})$. In section {\bf 2.}, we explain how to get a direct correspondence for the irreducible representations in $\R^{\varphi _0}_F(\underline{G})$ by parameters $(s,u,\Xi )$ in a given complex reductive group. The last section {\bf 3.} is devoted to the parametrization of representations of affine Hecke algebras, taking into account the local Langlands correspondence. At the end, corollary {\bf 3.5}, we give the final decomposition result (which does not depend on a final written account of the local Langlands correspondence). There are three annexes {\bf A}, {\bf B} and {\bf C}. In annex {\bf A}, it is explained how results in \cite{H3} generalize to the full orthogonal group which is not connected. In the annex {\bf B}, we give an account of the notion of tensor product in the context of linear abelian categories and show that it applies to the categories that we are considering. Unitary groups are treated in annex {\bf C}, although the results are used progressively in the main body of the paper.

Remark that only those results of this paper which apply to non quasi-split inner forms of $G$ are conditional: for orthogonal groups, some generalization of \cite{A} to inner forms may be required for \cite{M1, M3}, but this is forthcoming, and, for unitary groups, the case of the non split inner form has not been treated in \cite{M2, M3}, but the result is expected to be true and the proof not to be a major problem. Remark that \cite{A} is based on the stabilization of the twisted trace formula (or at least results which take part of this stabilization), but this has been accomplished recently thanks to a long series of papers by J.-L. Waldspurger.

One may expect that a similar pattern holds for a general $p$-adic reductive group. The method of Arthur based on the trace formula does not apply to the general case, but there are techniques on the Hecke algebra side developed by E. Opdam \cite{O1, O2} which may be a guide, once the results in \cite{H3} generalized to an arbitrary reductive $p$-adic group, to get a conjectural description of local $L$-packets for a general reductive $p$-adic group. A generalization of \cite{L2, L3} to groups which are not adjoint would also be helpful.

The author thanks S. Riche for many helpful discussions on Lusztig's work \cite{L2, L3}.

\null {\bf 1.} We will denote by ${^LG}$ the "$L$-group" of $G$, which means that its connected component is the Langlands dual group of the connected component of $G$ and that it is either a symplectic or a full (disconnected) orthogonal group. We will write $Z_{{^LG}}$ for the center of the connected component of ${^LG}$ (which is trivial if and only if $G$ is symplectic and of order two otherwise) and denote by $\iota :{^LG}\rightarrow GL_N(\Bbb C)$ the canonical embedding, i.e. $N$ equals $2n$ if $G$ is orthogonal and $N$ equals $2n+1$ if $G$ is symplectic. If $l$ is an integer between $1$ and $n$, $H_l$ will denote (the group of $F$ rational points) of a split classical group of semi-simple rank $l$ of the same type (symplectic, even or odd orthogonal) than $G$. The symbols $H_l^+$ and $H_l^-$ will have the appropriate meaning. We will also denote by $\iota $  the canonical embedding ${^LH}_l\rightarrow GL_L(\Bbb C)$, hoping that this will not be a source of confusion.

The connected component of $^LG$, the Langlands dual group, will be denoted $\widehat{G}$.

Let $W_F$ be the Weil group of $F$. It's the semi-direct product of the inertial subgroup $I_F$ with the cyclic subgroup generated by a Frob\'enius automorphism $Fr$, $W_F=\langle Fr\rangle\ltimes I_F$. A character $W_F$ is called \it unramified, \rm if it is trivial on $I_F$. By local class field theory, such a character is identified with a character of $F^{\times }$, trivial on the units of its ring of integers $O_F$, the character sending $Fr^{-1}$ to $q$ being identified with the absolute value $\vert\cdot\vert _F$.

We will call Langlands parameter for $G$ a continuous homomorphism $\rho $ of $W_F$ into ${^LG}$ which sends $Fr$ to a semi-simple element and assume in addition in the even-orthogonal case that the kernel of $det\circ\rho $ equals the Weil group of the splitting field of $G$. (It follows from the continuity that the image of $I_F$ is finite.) In the case, where $^LG$ is the odd orthogonal group, parameters with $det(\rho )\ne 1$ just correspond to representations of an identical copy of the symplectic group. A homomorphism $\rho :W_F\times SL_2(\Bbb C)\rightarrow {^LG}$ will be called a Langlands-Deligne parameter, if its restriction to the first factor is a Langlands parameter and the restriction to the second factor a morphism of algebraic groups. A Langlands or Langlands-Deligne parameter for ${^LG}$ will be called \it discrete, \rm if its image is not included in a proper Levi subgroup. Two Langlands or Langlands-Deligne parameters are said \it equivalent \rm if they are conjugated by an element of $^LG$. (Usually, one considers only conjugation by an element of $\widehat{G}$, but, as we take here for $G$ the full (non-connected) even orthogonal group, one has to take conjugation in $^LG$. For the other groups, this does not matter \cite{GGP, 8.1 (ii)}.)

If $\rho $ is an irreducible representation of $W_F$, the set of equivalence classes of representations of the form $\rho ^s:=\rho\vert\cdot\vert _F^s$, $s\in\Bbb C$, will be called the \it inertial class \rm of $\rho $. The group of unramified characters of $W_F$ acts on the inertial class of $\rho $ by torsion. We will denote by $t_{\rho }$ the order of the stabilizer of the equivalence class of $\rho $. If $\rho $ and $\rho '$ are in the same inertial class, then $t_{\rho }=t_{\rho '}$, and the definition of $t_{\rho }$ does not depend neither on the choice of $Fr$.

If $\rho $ is a self-dual representation, we will say that it is of type ${^LG}$, if it factors through a group of type ${^LG}$ (meaning that the image of $\rho $ is contained in an orthogonal group if ${^LG}$ is orthogonal and in a symplectic group if ${^LG}$ is symplectic). Otherwise, we will say that $\rho $ is not of type ${^LG}$. We stretch that the use of either of these notions will presume that $\rho $ is self-dual.

If $a$ is an integer $\geq 1$, $sp(a)$ will denote the unique irreducible representation of $SL_2(\Bbb C)$ of dimension $a$.

\null{\bf 1.1 Theorem:} \cite{M1, 1.5.1} \it 1) A Langlands-Deligne parameter $\varphi : W_F\times SL_2(\Bbb C)$ $\rightarrow {^LG}$ corresponds to a supercuspidal
representation of $G^+$ or $G^-$, if and only if $$\iota\circ\varphi =\bigoplus _{\rho\ not\ of\ type\ {^LG}}(\bigoplus _{k=1}^{a_{\rho }}(\rho\otimes sp(2k)))\oplus \bigoplus _{\rho\ of\ type\ {^LG}} (\bigoplus _{k=1}^{a_{\rho }}(\rho\otimes sp(2k-1))),$$

where the $\rho $ are irreducible representations and the $a_{\rho }$ non-negative integers.

2) Given $\varphi $ as in 1), denote by $z_{\varphi ,\rho ,k}$ the diagonal matrix in ${^LG}$ that acts by $-1$ on the space of the direct summand
$\rho\otimes sp(2k)$ (resp. $\rho\otimes sp(2k-1)$) of $\iota\circ\varphi $ and by $1$ elsewhere. Put $S_{\varphi }=C_{{^LG}}(Im(\varphi ))/C_{{^LG}}(Im(\varphi ))^{\circ }$. The elements
$z_{\varphi ,\rho ,k}$ lie in $C_{{^LG}}(Im(\varphi ))$ and their images $\ol{z}_{\varphi ,\rho ,k}$ generate the commutative group $S_{\varphi }$.

A pair $(\varphi ,\epsilon )$ formed by a discrete Langlands-Deligne parameter as in 1) and a character $\epsilon $ of $S_{\varphi }$ corresponds to a supercuspidal
representation of either $G^+$ or $G^-$, if and only if $\epsilon $ is alternating, i.e. $\epsilon (\ol{z}_{\varphi ,\rho ,k})=(-1)^{k-1}\epsilon (\ol{z}_{\varphi ,\rho ,1})$ with $\epsilon (\ol{z}_{\varphi ,\rho ,1})=-1$ for $\rho $ not of type ${^LG}$ and $\epsilon (\ol{z}_{\varphi ,\rho ,1})\in\{1, -1\}$ for $\rho $ of type ${^LG}$. It corresponds to a supercuspidal representation of $G^+$ if $\epsilon _{\vert Z_{{^LG}}}=1$ and to a supercuspidal representation of $G^-$ otherwise.

3) Suppose  that $\varphi $ satisfies the property in 1). Let $t_o$ be the number of $\rho $ of type ${^LG}$ with $a_{\rho }$ odd, put $t_o=1$ if there are none of them, and let $t_e$ be the number of the remaining $\rho $ of  type ${^LG}$ for which $a_{\rho }$ is even.

If $G$ is symplectic, there are $2^{t_o-1}2^{t_e}$ non isomorphic supercuspidal representations of $G^+$ with Langlands-Deligne parameter $\varphi $.

If $G$ is orthogonal, put $\epsilon _{\varphi ,\rho }=(-1)^{a_{\rho }(a_{\rho }+1)\over 2}$, if  $\rho $  is not of type ${^LG}$, and put $\epsilon _{\varphi ,\rho }=(-1)^{a_{\rho }\over 2}$, if  $\rho $  is of type ${^LG}$ and $a_{\rho }$ even. There exists a supercuspidal representation of $G^+$ with Langlands-Deligne parameter $\varphi $ if and only if either there is a $\rho $ of type ${^LG}$ with $a_{\rho }$ odd or  $\prod_{\rho } \epsilon_{\varphi ,\rho }=1$.

If the above existence condition is satisfied, the number of supercuspidal representations with Langlands-Deligne parameter $\varphi $ equals
$2^{t_o-1}2^{t_e}$ and all these representations of $G^+$ are non isomorphic.

The remaining alternating characters correspond to representations of $G^-$, remarking that there are $2^{t_0+t_e}$ alternating characters for orthogonal $G$.

\null Proof: \rm 1) and 2) are stated as this in the paper of Moeglin. Concerning 3), if $G$ is an orthogonal group, the theorem in the paper of Moeglin says that there is a supercuspidal representation of $G^+$ associated to $\varphi $, if and only if there exists an alternating character $\epsilon _{\varphi }$ corresponding to $\varphi $ which takes value $1$ on $-1$. The number of non isomorphic supercuspidal representations corresponding to $\varphi $ equals the number of alternating characters with this property.

For $\rho $ not of type ${^LG}$, there is a unique choice of an alternating character and its value on $-1$ is $\prod _{k=1}^{a_{\rho }} (-1)^k$. For $\rho $ of type ${^LG}$, there are always two choices of an alternating character. If $a_{\rho }$ is even and not divisible by 4 the value taken on $-1$ is always $-1$. If $a_{\rho }$ is divisible by $4$, the value taken on $-1$ is always $1$. If $a_{\rho }$ is odd, there is one alternating character which takes value $1$ on $-1$ and another one which takes value $-1$ on $-1$.

One concludes by remarking that, if there are alternating characters attached to a $\rho $ which take respectively value $1$ and $-1$ on $-1$, then one can of course always find an alternating character for $\varphi $ with value $1$ on $-1$.

If $G$ is symplectic, one can conclude as above, after having observed that there is then always a $\rho $ of type ${^LG}$ with $a_{\rho }$ odd. \hfill{\fin 2}

\null {\bf 1.2.} \it Definition: \rm We will fix for the rest of the paper in each inertial class $\o $ of an irreducible representation of $W_F$ a base point $\rho_{\so }$. It will always be assumed to be the equivalence class of a unitary representation, which is in addition self-dual if $\o $ contains such an element. In this last case, we take $\rho $ of the same type as ${^LG}$ if there is such a representation in $\o $. This base point will be called in the sequel a \it normed \rm representation (w.r.t. ${^LG}$).

A discrete Langlands parameter $\tau :W_F\rightarrow {^LG}$ will be called \it normed, \rm if $\iota\circ\tau $ is the direct sum of inequivalent normed
representations of $W_F$. If $^LM\simeq GL_{k_1}(\Bbb C)\times\cdots\times GL_{k_r}(\Bbb C)\times {^LH}_l$ is a Levi subgroup of ${^LG}$, then a discrete Langlands parameter $\varphi :W_F\rightarrow{^LM}$ is called normed if it is of the form $\gamma\mapsto (\rho_1(\gamma ), \dots, \rho _r(\gamma ),\tau (\gamma ))$, where the $\rho _i$ are irreducible normed representations and $\tau $ is a discrete normed Langlands parameter. A Langlands parameter $\varphi :W_F\rightarrow {^LG}$ will be called normed, if there is a minimal Levi subgroup $^LM$ of ${^LG}$ containing the image of $\varphi $, such that $\varphi $ is a discrete normed Langlands parameter with respect to $^LM$.

If $s$ is a semi-simple element  in $GL_N(\Bbb C)$, then $\chi _s$ will denote the unramified character of $W_F$, such that $\chi _s(Fr)=s$. If $\varphi $ is a normed Langlands parameter and $s$ is a semi-simple element in $C_{GL_N(\Bbb C)}(Im(\varphi ))$ such that $(\iota\circ\varphi)\chi _s$ is a self-dual representation of $W_F$ of the same type than $^LG$, then we will denote by $\varphi _s$ the Langlands parameter $W_F\rightarrow{^LG}$ such that $\iota\circ\varphi_s$ is equivalent to the representation $\gamma\mapsto\varphi _s(\gamma )=\varphi(\gamma )\chi_s(\gamma )$ of $W_F$.

The set of equivalence classes of Langlands parameters of the form $\varphi _s$ with $s$ as above will be called the \it inertial orbit \rm of $\varphi $.

If $\rho '$ is an irreducible representation of $W_F$, $m(\rho ';\varphi )$ will denote the multiplicity of $\rho '$ (up to equivalence) in the representation
$\rho _1\oplus\cdots\oplus\rho_r\oplus\rho_r^{\vee }\oplus\cdots\oplus\rho_1^{\vee }\oplus(\iota\circ\tau )$.

\null \it Remark: \rm Normed Langlands parameters will play a similar role than the trivial Langlands parameter for the set of unipotent representations. Unfortunately, it seems not possible to fix in general a "canonical" base point (see also remark after {\bf 1.4}). The choice of a base point has no influence on the essentially intrinsic nature of our results.

\null {\bf 1.3 Proposition:} \it Two Langlands parameters $W_F\rightarrow{^LG}$ are in the same inertial orbit, if and only if their restrictions to the inertial subgroup $I_F$ are conjugate by an element of $^LG$.

\null Proof: \rm For the proof, it is enough to consider the case where one of the two Langlands parameters is normed. So, let $\varphi $ and $\varphi '$ be two Langlands parameters.

Suppose first that $\varphi $ and $\varphi '$ are in the same inertial orbit as defined above. Then, as representations of $I_F$, $\iota\circ\varphi _{\vert I_F}$ and $\iota\circ\varphi' _{\vert I_F}$ are certainly isomorphic. One deduces from this, analog to \cite{GGP, 8.1 (ii)} in the case of Langlands parameters, that $\varphi _{\vert I_F}$ and $\varphi '_{\vert I_F}$ are conjugate by an element of $^LG$. This proves one direction.

For the other direction suppose that $\varphi _{\vert I_F}$ and $\varphi '_{\vert I_F}$ are conjugate by an element of $^LG$. Then, as representations of $I_F$, $\iota\circ\varphi _{\vert I_F}$ and $\iota\circ\varphi' _{\vert I_F}$ are certainly isomorphic. Remark that an irreducible representation of $W_F$ is determined, up to twist by an unramified character, by its restriction to $I_F$ and more particular by an irreducible component of this restriction. It follows that for each irreducible component $\rho $ of $\iota\circ\varphi $, there is an irreducible component $\rho '$ of $\iota\circ\varphi '$ such that the two irreducible representations have a common irreducible component when restricted to $I_F$. In addition, $\rho '$ is in the inertial class of $\rho $ (as irreducible representation of $W_F$) and the restrictions $\rho _{\vert I_F}$ and $\rho '_{\vert I_F}$ are isomorphic. The same is true, if one starts from an irreducible component of $\iota\circ\varphi '$. This induces a bijection between irreducible components of $\iota\circ\varphi $ and irreducible components of $\iota\circ\varphi '$, sending an irreducible component of $\iota\circ\varphi $ to an irreducible representation in its inertial orbit. It follows that there is a semi-simple element $s$ in $C_{GL_N(\Bbb C)}(Im(\varphi ))$ such that $\iota\circ\varphi '$ is isomorphic to the representation $\gamma\mapsto\varphi(\gamma )\chi_s(\gamma )$. By definition, this means that $\varphi '$ is in the inertial orbit of the normed Langlands parameter $\varphi $.\hfill{\fin 2}

\null {\bf 1.4 Proposition:} \it Let $\o $ be the inertial orbit of an irreducible representation of $W_F$ and denote by $\rho _{\so }$ the normed representation in its inertial orbit. The map $\o\rightarrow\Bbb C$, defined by $\rho \mapsto f_{\rho }:=\rho (Fr^{t_{\rho }})\rho _{\so }(Fr^{t_{\rho }})^{-1}$, where $\rho $ denotes a representative of the equivalence class which is an unramified twist of $\rho _{\so }$, is a well-defined bijection. If $\rho _{\so }$ is self-dual, then $\rho $ is self-dual, if and only if $f_{\rho }\in\{\pm 1 \}$.

\null Proof: \rm By definition, there is a complex number $s$ such that $\rho =\rho _{\so }\otimes\vert\cdot\vert_F^s$. It follows from this that the map $\rho\mapsto f_{\rho }$ is well defined and that, for $\rho $ as above, $f_{\rho }=q^{-st_{\rho }}$. One sees that the map is surjective. Put $\rho '=\rho _{\so }\otimes\vert\cdot\vert _F^{s'}$. Then $f_{\rho }=f_{\rho '}$ is equivalent to $q^{(s-s')t_{\rho }}=1$. This implies that $\vert\cdot\vert_F^{s'-s}$ stabilizes the equivalence class of $\rho $. Consequently, $\rho'=\rho\otimes\vert\cdot\vert_F^{s'-s}\simeq\rho $. So, the map is also injective.

Assume now $\rho _{\so }$ self-dual and that $\chi $ is an unramified character such that $\rho:=\rho_{\so }\otimes\chi $ is also self-dual. This implies that $\chi ^2$ stabilizes $\rho _{\so }$ and consequently one has $f_{\rho }^2=1$. One sees that it is enough to twist $\rho _{\so }$ by $\vert\cdot\vert_F^{i\pi \over t_{\rho }log(q)}$ to get a representation $\rho '$ that satisfies $f_{\rho '}=-1$.\hfill{\fin 2}

\null \it Remark: \rm In general, it is not possible to distinguish an element $\rho $ of $\o $, such that $f_{\rho }=1$, even if $\rho $ is self-dual: although
$\rho $ is then induced from an irreducible representation of the Weil group of an unramified extension of degree $t_{\rho }$, there is no reason why $\rho (Fr^{t_{\rho }})$ should be a scalar. (For example, in the case $t_{\rho }=1$, conjugation by $Fr $ gives an isomorphic representation of the restriction to the inertia group, but it is not necessarily the same representation.) That is the reason, why we had to make a choice in our definition of a normed representation.

\null {\bf 1.5} \it Definition: \rm If $\rho $ is self-dual, we will denote by $\rho _-$ the unique element in its inertial orbit such that $f_{\rho _-}=-1$. When $H$ is an even orthogonal quasi-split group and $\zeta $ is the quadratic character of $W_F$ whose kernel corresponds to the splitting field of $H$, we will denote by $H_-$ the quasi-split outer form of $H$ whose splitting field corresponds to the kernel of $\zeta _-$. We will leave $H_-$ undefined when $H$ is odd-orthogonal and put $H_-=H$ if $H$ is symplectic, although $H_-$ will be distinguished from $H$. We will also write $H_+$ for $H$. The notation $H_\pm^\pm$ will then have the appropriate meaning.

\null\it Remark: \rm If $\rho $ is self-dual, $\rho $ is either orthogonal or symplectic. However, it happens that $\rho _-$ is not of the same type (i.e. orthogonal or
symplectic) than $\rho $, and both cases, $\det(\rho )=\det(\rho _-)$ and $\det(\rho )\ne\det(\rho _-)$, happen. (Examples can be easily deduced from \cite{Mo, theorem 1}.) By our convention of a normed representation, this can only appear if $\rho _0$ is of type ${^LG}$.

\null{\bf 1.6} If $\varphi _0:W_F\rightarrow {^LG}$ is a normed Langlands parameter, denote by $supp(\varphi _0)$ the set of irreducible representations $\rho $ of
$W_F$, up to isomorphism, with $m(\rho ;\varphi _0)\ne 0$ and by $supp'(\varphi _0)$ the subset formed by those representations which are selfdual. We will put an equivalence $\sim $ on $supp(\varphi _0)$ defined by $\rho\sim\rho^{\vee }$.

Denote by $\hbox{\main S}(\varphi _0)$ the set of families of pairs $(a_{\rho ,+}, a_{\rho ,-})_{\rho }$ of non-negative integers indexed by $supp'(\varphi _0)$, such that
$$m(\rho;\varphi _0)\geq \cases a_{\rho,+}(a_{\rho,+}+1)+a_{\rho,-}(a_{\rho,-}+1), & \hbox{\rm if $\rho$ and $\rho _-$ not of type ${^LG}$},\cr
a_{\rho,+}^2+a_{\rho,-}(a_{\rho,-}+1), & \hbox{\rm if $\rho$ of type ${^LG}$, but not $\rho _-$},\cr
a_{\rho,+}^2+a_{\rho,-}^2, & \hbox{\rm if $\rho$ and $\rho _-$ of type ${^LG}$},\cr\endcases$$
with the additional condition that the terms of both sides in the above inequalities have same parity (if $\rho $ and $\rho _-$ are both not of type ${^LG}$, this is always satisfied).

Put $\kappa '_{\rho }=1$ if $\rho $ is of type ${^LG}$ and $\kappa '_{\rho }=0$ otherwise. If $S=(a_{\rho ,+}, a_{\rho ,-})_{\rho }$ lies in $\hbox{\main S}(\varphi _0)$, then the dimension $L_S$ of the representation
$$\oplus _{\rho\in supp'(\varphi _0)}(\oplus _{k=1}^{a_{\rho,+}}(\rho\otimes sp(2k-\kappa'_{\rho }))\oplus  (\oplus _{k=1}^{a_{\rho ,-}}\rho_-\otimes sp(2k-\kappa'_{\rho _-})))$$ has the same parity than $N$. If we denote by $l_S$ the semi-simple rank of the group $H_{l_S}$ whose Langlands dual embeds canonically into $GL_{L_S}(\Bbb C)$, then there is, up to equivalence, a unique discrete Langlands-Deligne parameter $\varphi ^S:W_F\rightarrow\ ^L{H_{l_S}}$ \cite{GGP, 8.1.ii} (as we consider the full orthogonal group, the restriction for the even orthogonal group does not apply), such that the above representation is equivalent to $\iota\circ\varphi ^S$. Denote by $\widehat{S}$ the set of alternating characters of $S_{\varphi ^S}$ (see theorem {\bf 1.1}, 2) for the definition of this group) and, for $\epsilon\in\widehat{S}$, by $\epsilon _Z$ its restriction to $Z_{{^LG}}$ (which can be $1$ or $-1$). Put $d_S=+$ or $d_S=-1$ according to whether $det(\varphi ^S)=det(\varphi_0)$ or not. (Remark that in the latter case necessarily $det(\varphi ^S)=det(\varphi_0)_-$ seen as representation of $W_F$.) Then $\varphi ^S$ defines a Langlands-Deligne parameter for the quasi-split group $H_{l_S, d_S}$. Write $\tau _{S,\epsilon }$ for the irreducible supercuspidal representation of $H_{l_S, d _S}^{\epsilon _Z}$ which corresponds to the Langlands-Deligne parameter $\varphi ^S$ and the alternating character $\epsilon $ of $S_{\varphi ^S}$.

Denote by $k_{\rho }$ the dimension of $\rho $ and put $m_\pm(\rho;\varphi ^S)=m(\rho;\varphi ^S)+m(\rho_-;\varphi ^S)$. Define $^LM_S$ to be the Levi subgroup of ${^LG}$ that is isomorphic to
$$\eqalign{\prod _{\rho\in (supp(\varphi _0)-supp'(\varphi _0))/\sim }&GL_{k_{\rho }}(\Bbb C)^{m(\rho;\varphi _0)}
\times\cr &\times\prod _{\rho\in supp'(\varphi _0)}GL_{k_{\rho }}(\Bbb C)^{[(m(\rho;\varphi _0)-m_\pm(\rho;\varphi ^S))]/2}\times {^LH}_{l_S}.\cr}$$
(Here $/\sim $ stands for the equivalence classes w.r.t. the relation $\rho\sim\rho ^{\vee }$ defined above.) Let $\varphi _S$ be the discrete Langlands-Deligne parameter $W_F\times SL_2(\Bbb C)\rightarrow {^LM}_S$
such that $$\iota\circ\varphi _S=\oplus _{\rho\in supp(\varphi _0)}{[m(\rho;\varphi _0)-m(\rho;\varphi ^S)]}\rho\oplus(\iota\circ\varphi ^S).$$

Denote by $M_S^{\epsilon _Z}$ the standard Levi subgroup of $G_{d_S}^{\epsilon _Z}$ with $L$-group $^LM$. It is isomorphic to a product of general linear groups with one factor isomorphic to $H_{l_S, d _S}^{\epsilon _Z}$.   For $S\in\hbox{\main S}(\varphi _0),\ \epsilon\in\widehat {S}$, denote by $\sigma _{S,\epsilon }$ the supercuspidal representation of $M_S^{\epsilon _Z}$ which corresponds to $\varphi _S$ and $\epsilon $ (i.e. the factor $H_{l_S, d_S}^{\epsilon _Z}$ acts by $\tau _{S,\epsilon }$) and by $\o _{S,\epsilon }$ the corresponding inertial orbit, i.e. $\o _{S,\epsilon }$ is the set of equivalence classes of representations of $M_S^{\epsilon _Z}$ which are unramified twists of $\sigma _{S,\epsilon }$.

In general, if $\sigma '$ is an irreducible supercuspidal representation of $M'$, we will denote by $\varphi _{\sigma '}$ the Langlands-Deligne parameter of $\sigma '$ obtained by applying {\bf 1.1} and the local Langlands correspondence to the $GL_k$.

\null{\bf Theorem:} \it The family $(M_S^{\epsilon _Z},\o _{\sigma _{S,\epsilon }})_{S,\epsilon }$ exhausts the set of inertial orbits of supercuspidal representations $\sigma '$ of
standard Levi subgroups $M'$ of $G$, $G_{-}$, $G^{-}$ and $G_{-}^{-}$ with $\widehat{M'}\supseteq \widehat{M}$, such that  $(\varphi _{\sigma '})_{\vert W_F}$ lies in the inertial orbit of $\varphi _0$.

One has $(M_S,\o _{\sigma _{S,\epsilon }})=(M_{S'},\o _{\sigma _{S',\epsilon '}})$, if and only if $(S,\epsilon )=(S',\epsilon ')$.

\null
\it Proof: \rm The first part follows directly from the constructions and theorem {\bf 1.1}. For the second part: to have $(M_S,\o _{\sigma _{S,\epsilon }})=(M_{S'},\o _{\sigma _{S',\epsilon '}})$, one needs $l_S=l_{S'}$ and $\tau _{S,\epsilon }=\tau _{S',\epsilon '}$, but then the other factors of $\sigma _{S,\epsilon }$ and $\sigma _{S',\epsilon '}$ must be unramified twists of each other. \hfill{\fin 2}

\null{\bf 1.7} We summarize below the (partly expected) properties of the local Langlands correspondence for $G$, which is in quite good shape now (see remarks below). If $\o $ is the inertial orbit of a supercuspidal representation of a Levi subgroup of $G_\pm^\pm $, we will denote by $Rep_{\so }(G_\pm^\pm )$ the corresponding Bernstein component of $Rep(G_\pm^{\pm })$ \cite{BD}.

If $\varphi _0:W_F\rightarrow {^LG}$ is a normed Langlands parameter for $G$, put $\S_+(\varphi _0)=\{S\in\S(\varphi _0)\vert\ det(\varphi _S)=det(\varphi _0)\}$, $\S_-(\varphi _0)=\{S\in\S(\varphi _0)\vert \ det(\varphi _S)\ne det(\varphi _0)\}$, and, for $S\in\hbox{\main S}(\varphi _0)$, $\widehat{S}^\pm=\{\epsilon\in\widehat{S}\vert \epsilon _Z=\pm 1\}$,   $$\R _{F,\pm }^{\varphi _0, \pm}(\underline{G})=\sum _{S\in \hbox{\smain S}_{\pm}(\varphi _0), \epsilon\in\widehat{S}^{\pm }}Rep_{\so_{S,\epsilon }}(G_\pm^\pm).$$ Denote by $\R _F^{\varphi _0,\cdot }(\underline{G})$ the direct sum of $\R _{F,+}^{\varphi _0,\cdot }(\underline{G})$ and $\R _{F,-}^{\varphi _0,\cdot }(\underline{G})$ and by $\R _F^{\varphi _0}(\underline{G})$ the direct sum of $\R _{F}^{\varphi _0,+}(\underline{G})$ and $\R _{F}^{\varphi _0,-}(\underline{G})$.

\null {\bf Local Langlands Correspondence.}

\it For a fixed normed Langlands parameter $\varphi _0:W_F\rightarrow {^LG}$ for $G$, the set of equivalence classes of pairs $(\varphi, \Xi)$ with $\varphi :W_F\times SL_2(\Bbb C)\rightarrow\ {^LG}$ a Deligne-Langlands parameter with $\varphi _{\vert W_F}$ in the inertial orbit of  $\varphi _0$, and $\Xi$ an irreducible representation of $C_{\widehat{G}}(Im(\varphi ))/(C_{\widehat{G}}(Im(\varphi )))^0$, is in natural bijection with $\R _F^{\varphi _0}(\underline{G})$ .

Pairs $(\varphi, \Xi)$ with $\Xi _{\vert Z_{{^LG}}}=1$ (resp. $\Xi _{\vert Z_{{^LG}}}=-1$) correspond to representations of $G_\pm^+$ (resp. $G_\pm^-$), the ones with $det(\varphi )=det(\varphi _0)$ (resp. $det(\varphi )\ne det(\varphi _0))$ to representations of $G_+^\pm$ (resp. $G_-^\pm$), those with $\varphi $ discrete to square integrable representations and the ones with $\varphi (W_F)$ bounded to tempered representations.

All smooth irreducible representations of $G^+$ and $G^-$ appear when $\varphi _0$ varies.

In addition, the following  equalities of local constants hold: if $^LM$ is the standard Levi subgroup  of a maximal standard parabolic subgroup $^LP$ of
${^LG}$, denote by $r_1$, $r_2$ the irreducible components of the regular representation of $^LM$ on the Lie algebra of the unipotent radical of $^LP$. Let $\pi $ be an irreducible smooth representation of the corresponding maximal Levi subgroup $M$ of $G$ and $\varphi _{\pi }:W_F\times SL_2(\Bbb C)\rightarrow{^LM}$ its Langlands-Deligne parameter. Then, the local factors defined by the Langlands-Shahidi method satisfy, for $i=1,2$,
$$\gamma (r_i\circ\varphi _{\pi },s)=\gamma(\pi ,r_i,s),\ \epsilon (r_i\circ\varphi _{\pi },s)=\epsilon(\pi , r_i,s)\ \hbox{\it and}\ L(r_i\circ\varphi _{\pi },s)=L(\pi, r_i,s).$$

\null Remark: \rm (i) The representations $\Xi $ have to be taken relative to the group of connected components defined by the centralizer of $Im(\varphi )$ in $\widehat{G}$, although the image of $\varphi $ may not lie in $\widehat{G}$. Remark that this difference did not matter for {\bf 1.1}.

(ii) It is explained in \cite{Sh, 8} how to define the local factors for non generic representations and also for representations of inner forms (see also \cite{H2, section 4}).

(iii) If $G$ is quasi-split and $^LM=GL_k(\Bbb C)\times {^LH}_l$, then $r_1$ is the standard representation $id_{GL_k(\Bbb C)}\otimes\iota$ and $r_2=Sym^2\circ id_{GL_k(\Bbb C)}$ or $\wedge^2\circ id_{GL_k(\Bbb C)}$, depending if ${^LG}$ is symplectic or orthogonal.

(iv) The local Langlands correspondence for classical groups is in quite good shape now in consequence of the work of J. Arthur who describes the discrete series and
tempered representations with their $R$-groups (see \cite{A} for the split case, the case of inner forms is forthcoming) and the work of C. Moeglin \cite{M1, M3}. Results on the preservation of local factors for symplectic and orthogonal Galois representations have been established by Cogdell-Shahidi-Tsai \cite{CST}, completing work of G. Henniart \cite{He2}. However, no final account has been written on all this yet.

Remark that, for $G=SO_{2n+1}(F)$, the case of $\varphi _0=1$ has been solved in \cite{L2} with additional work in \cite{W}. As the constructions in \cite{W} are compatible with \cite{H2} and the local Langlands correspondence for quasi-split tori is known to preserve local factors, it follows from \cite{H2} that local factors are preserved for this correspondence (which coincides, at least on the level of Langlands-Deligne parameters, with the one in {\bf 1.7}).

(v) For the  group $Sp_4(F)$, W.-T. Gan and S. Takeda gave in \cite{GT} properties for the local Langlands correspondence, which makes it unique.

\null{\bf 1.8} The following result of \cite{H4} is obtained by linking the results of C. Moeglin to \cite{H3} (see the remark after theorem {\bf A.7} in the annex for the case of the non connected orthogonal group). Recall that it can well happen that $\rho $ is orthogonal and $\rho _-$ symplectic or vice versa \cite{Mo}. The terminology for affine Hecke algebras with parameters used below is the one from \cite{L1}, after evaluation in $q^{1/2}$ as done in \cite{L2, L3}.

Recall the equivalence relation on $supp(\varphi _0)$ given by $\rho\sim\rho^{\vee }$ introduced in {\bf 1.6}.

\null {\bf Theorem:} \cite{H4, Theorem 5.2 and remark thereafter} \it  Let $\varphi _0$ be a normed Langlands parameter, $S\in \hbox{\main S}(\varphi _0)$,
$S=(a_{\rho ,+}, a_{\rho ,-})_{\rho }$  and $\epsilon\in\widehat{S}$.

The category $Rep_{\so _{S,\epsilon }}(G_{d_S}^{\epsilon _Z})$ is equivalent to the category of right modules over the tensor product
$\otimes _{\rho\in supp(\varphi _0)/\sim }\H _{\varphi _0,S,\rho }$ where $\H _{\varphi _0,S,\rho }$ are extended affine Hecke algebras of the following type:

- if $\rho $ is not self-dual, $\H_{\varphi _0,S,\rho }$ is an affine Hecke algebra with root datum equal to the one of $GL_{m(\rho ;\varphi _0)}$ and equal parameters $q^{t_{\rho }}$;

- if $\rho $ and $\rho _-$ are both of the same type than ${^LG}$ and $a_{\rho ,+}=a_{\rho ,-}=0$, then $\H_{\varphi _0,S,\rho  }$ is the semi-direct product of an affine Hecke algebra with root datum equal to the one of $SO_{m(\rho ;\varphi _0)}$ and equal parameters $q^{t_{\rho }}$ by the group algebra of a finite cyclic group of order $2$, which acts by the outer automorphism of the root system;

- otherwise, putting $\kappa _{\rho ,\pm 1}=0$ (resp. $=1$) if  $\rho _\pm$ is of type ${^LG}$ (resp. not of type ${^LG}$), $\H_{\varphi _0,S,\rho  }$ is an affine Hecke algebra with root datum equal to the one of $SO_{m(\rho ;\varphi _0)-m_\pm(\rho  ;\varphi ^S)+1}$ and unequal parameters
$q^{t_{\rho }},\dots ,q^{t_{\rho }}, q^{t_{\rho } (a_{\rho ,+}+a_{\rho ,-} +({{\kappa _{\rho ,+}+\kappa _{\rho ,-}}\over 2}))};$
$q^{t_{\rho } \vert a_{\rho ,+}-a_{\rho ,-} +({{\kappa _{\rho ,+}-\kappa _{\rho ,-}}\over 2})\vert}$, remarking that $m(\rho ;\varphi _0)-m_\pm(\rho  ;\varphi ^S)+1$ is necessarily an odd number. \rm

\null\it Remark: \rm (i) If $\rho $ is not of type ${^LG}$ and $a_{\rho ,+}=a_{\rho ,-}=0$, then it is well known that the affine Hecke algebra $\H_{\rho }$ expressed above
is isomorphic to an affine Hecke algebra with root datum equal to the one of $Sp_{m(\rho ;\varphi _0)}$ and equal parameter $q^{t_{\rho }}$.

(ii) The $\kappa _{\rho ,\cdot }$ in the above theorem is related to the $\kappa' _{\rho }$ in {\bf 1.6} by the relation $\kappa '_{\rho }=1-\kappa _{\rho ,+}$ and $\kappa '_{\rho _-}=1-\kappa _{\rho ,-}$.

(iii) By \cite{H5}, the above equivalence of categories preserves temperedness and discreteness (in the definition for square integrability modulo the "center" for Hecke algebra-representations, the "center"  has of course to be taken trivial, so that there are no discrete series representations if the based root system which defines the Hecke algebra has a factor of type $A_n$). Unitarity is conjectured.

\null {\bf 1.9} The notion of a tensor product of linear abelian categories is treated in \cite{D} and recalled in the annex {\bf B}. It applies to the category of modules of finite presentation over a coherent $\Bbb C$-algebra and in particular to the category of finitely generated modules over a noetherian $\Bbb C$-algebra (for ex. an extended affine Hecke algebra or a finite tensor product of such algebras (cf. {\bf B.3})). Denote by  $\R^{\varphi _0}_F(\underline{G})_f$ the full subcategory of  $\R^{\varphi _0}_F(\underline{G})$ whose objects are the finitely generated representations (by \cite{BD, 3.10}, these are precisely the representations which are admissible relative to the action of the Bernstein center) and  recall the equivalence relation on $supp(\varphi _0)$ given by $\rho\sim\rho^{\vee }$ introduced in {\bf 1.6}.

\null {\bf Corollary:} \it The category $\R^{\varphi _0}_F(\underline{G})_f$ is equivalent to the category $$\bigoplus_{S\in \sS(\varphi _0), \epsilon\in \widehat{S}^\pm} (\bigotimes _{\rho\in supp(\varphi _0)/\sim} (right-\H_{\varphi _0,S,\rho  }-modules)_f\ ).$$

\null Proof: \rm It follows from {\bf B.2} and {\bf B.3} that the tensor product exists and can be applied by {\bf 1.8} to the category $\R^{\varphi _0}_F(\underline{G})$ (defined in {\bf 1.7}) to give the statement of the corollary. \hfill{\fin 2}

\null{\bf 1.10} The above results generalize to pure inner forms of quasi-split unitary groups over $F$, as remarked in annex {\bf C}, {\bf C.1} - {\bf C.5}.

\null {\bf 2.} The object of this section is to relate,  for a given normed Langlands parameter $\varphi _0$, the Deligne-Langlands-Lusztig parameters for the irreducible representations in $\R _F^{\varphi _0}(G)$ to data given by semi-simple and unipotent elements in a given complex group. Recall that we have put in {\bf 1.6} an equivalence relation on the set of normed representations of $W_F$  (up to isomorphism) defined by $\rho \sim \rho ^{\vee }$.

\null
{\bf 2.1 Proposition:} (\cite{GGP, section 4})\it

Let ${^LM}=GL_{k_1}(\Bbb C)\times\cdots\times GL_{k_r}(\Bbb C)\times {^LH_l}$ be a  standard Levi subgroup of ${^LG}$, $\varphi :W_F\rightarrow {^LM}$ a discrete normed Langlands parameter, $\iota\circ\varphi =\rho_1\oplus\cdots\oplus\rho _r\oplus(\iota\circ\tau )$.

Then, one has $C_{{^LG}}(Im(\varphi ))=\prod_{\rho } H_{\rho ;\varphi }(m(\rho ;\varphi ))$, the product going over representatives of the equivalence classes of irreducible normed representations $\rho $ of $W_F$, while the $H_{\rho ;\varphi }(m)$ are complex classical groups with $H_{\rho ;\varphi }(m)$ isomorphic to $GL_m(\Bbb C)$ if $\rho $ is not self-dual, to $Sp_m(\Bbb C)$ if $\rho $ is not of type ${^LG}$ and to $O_m(\Bbb C)$ if $\rho $ is of type ${^LG}$ (with the convention $O_1(\Bbb C)=\{\pm 1\}$ if $m=1$).

Finally, $C_{GL_N(\Bbb C)}(Im(\varphi ))=\prod_{\rho } G_{\rho ;\varphi }(m(\rho ;\varphi ))$, the product going over representatives of the equivalence classes of irreducible representations $\rho $ of $W_F$, while $G_{\rho ;\varphi }(m)$ is isomorphic to $GL_m(\Bbb C)$, if $\rho $ is self-dual, to $GL_m(\Bbb C)\times GL _m(\Bbb C)$ if $\rho $ is not self-dual, and the group $G_{\rho ;\varphi }(m)$ contains $H_{\rho ;\varphi }(m)$ in each case.

On the other side,  $C_{GL_N(\Bbb C)}(Im(\varphi_s)) \subseteq \prod_{\rho }G_{\rho ;\varphi }(m(\rho ;\varphi ))$ for every unramified twist $\varphi _s$ of $\varphi $ with $s$ in the centralizer of $Im(\varphi )$. \rm

\null{\bf 2.2} Recall that the invariant $f_{\rho }$ of an irreducible representation $\rho $ of $W_F$ has been defined in {\bf 1.4}.

\null {\bf Lemma:} \it Let ${^LM}\simeq GL_{k_1}(\Bbb C)\times\cdots\times\GL_{k_d}(\Bbb C)\times {^LH_l}$ be a standard Levi subgroup of ${^LG}$ and let $\varphi :W_F\rightarrow{^LM}$ be a discrete Langlands parameter, $\gamma\mapsto (\rho _1(\gamma ),\dots ,\rho _d(\gamma ),\tau (\gamma ))$. Denote by $\varphi _0$ the normed Langlands parameter associated to $\varphi $. Write $\iota\circ\tau =\tau _1\oplus\cdots\oplus\tau _r$ for the decomposition of $\iota\circ\tau $ into irreducible representations. Denote by $s_{\tau }$ the element of $C_{^LH_l}(Im(\tau ))$ which corresponds to the diagonal matrix $(f_{\tau _1},\dots ,f_{\tau _r})$ and by $s_{\varphi }$ the element of $C_{{^LM}}(Im(\varphi ))$ that corresponds to $(f_{\rho _1},\dots , f_{\rho _d}, s_{\tau })$.

The element $s_{\varphi }$ lies in $C_{{^LG}}(Im(\varphi _0))$ and in $C_{{^LG}}(Im(\varphi ))$.

Suppose: if $\rho _i$ is self-dual, then it is of the same type than the normed representation in its inertial orbit. Then, $C_{{^LG}}(Im(\varphi ))$ and $C_{C_{{^LG}}(Im(\varphi _0))}(s_{\varphi })$ are canonically isomorphic.

\null
Remark: \rm As $\varphi $ and consequently $\tau $ are discrete, the representations $\tau _i$ are all of type ${^LG}$ and non isomorphic. In addition, $f_{\tau _i}\in\{\pm 1\}$ for $i=1,\dots ,r$.

\null \it Proof: \rm By the proposition {\bf 2.1}, one has $C_{{^LG}}(Im(\varphi _0))=\prod _iGL_{l_i}(\Bbb C)\times \prod_jSp_{m_j}(\Bbb C)$ $\times \prod_kO_{n_k}(\Bbb C)$, where the first product goes over the $\rho _i$ which are not self-dual, the second one over the $\rho _j$ which are not of the same type than ${^LG}$ and the third one over the $\rho _k$ which are of the same type than ${^LG}$ . The centralizer of $Im(\varphi _0)$ is determined by the partition of the summands of $\iota\circ\varphi _0$ obtained by putting together representations which are either isomorphic or isomorphic to the dual of the other one. The different parts of this partition of the summands  of $\iota\circ\varphi _0$ give then rise to factors which are respectively isomorphic to $GL_{l_i}(\Bbb C)$, $Sp_{m_j}(\Bbb C)$ or $O_{n_k}(\Bbb C)$ depending if the representations in the part are not self-dual, orthogonal or symplectic, where $l_i$ denotes half the number of elements in the corresponding part and $m_j$ and $n_k$ the total number of elements in the part. The analog result holds for the centralizer of $Im(\varphi )$. As $\varphi _0$ is normed, it is clear that the partition of $\iota\circ\varphi$ is finer than the partition of $\iota\circ\varphi _0$.

Writing $C_{{^LG}}(Im(\varphi _0))$ as above as a product, one sees that the centralizer of $s_{\varphi }$ in $C_{{^LG}}(Im(\varphi _0))$ is the product of the centralizers of the components of $s_{\varphi }$ in the different factors. So, one can reduce to consider the following three cases:

(i) All summands of  $\iota\circ\varphi $ are in the same inertial orbit and the normed representation in this orbit is not self-dual. In particular, $\tau $ is trivial.

(ii) All summands of $\iota\circ\varphi $ are in the same inertial orbit and the normed representation in this orbit is not of type $^LG$. In particular, $\tau $ is trivial.

(iii) All summands of $\iota\circ\varphi $ are in the same inertial orbit and the normed representation in this orbit is of type $^LG$. In particular, $\iota\circ\tau $ is either trivial or equal to an element of the inertial orbit of this normed representation of type $^LG$.

In all three cases, the centralizer of $s_{\varphi }$ is determined by the partition of the coefficients of $s_{\varphi }$, obtained by putting equal coefficients in the same part. By proposition {\bf 1.4}, equal coefficients correspond to equal summands of $\iota\circ\varphi $, so that the two partitions correspond canonically to each other and have the same number of elements. The factors of the centralizer of $s_{\varphi }$ which correspond to the different parts of the partition of the coefficients of $s_{\varphi }$ are all general linear groups of order equal to the length of the partition in case (i). In case (ii) and (iii) they are general linear groups if the coefficients are $\ne\pm 1$ and groups of the same type than the group in the other cases, as by our assumption the appearance of groups of another type is excluded. This proves the proposition. \hfill{\fin 2}

\null{\bf 2.3 Lemma:} \it With the same notations as in {\bf 2.2}, assume that $\rho $ is an irreducible representation of $W_F$ of type ${^LG}$, such that $\rho _-$ is not of type ${^LG}$ and  $\iota\circ\varphi _0\simeq m\rho $.

Then, $C_{GL_N(\Bbb C)}(Im(\varphi _0))$ is canonically isomorphic to $GL_m(\Bbb C)$, while $C_{{^LG}}$ $(Im(\varphi _0))$ is isomorphic to $O_m(\Bbb C)$.

Define $s_{\varphi }$ as in {\bf 2.2}.  The element $s_{\varphi }$ lies in $C_{{^LG}}(Im(\varphi _0))$ and in $C_{{^LG}}(Im(\varphi ))$.

Write $s_{\varphi }=diag(x_1,\dots ,x_{[{m\over 2}]},\widehat{1}, x_{[{m\over 2}]}^{-1},\dots , x_1^{-1})\in GL_m(\Bbb C)$ (with $1$ appearing only when $m$ is odd and $[{m\over 2}]$ denoting the integer part of ${m\over 2}$). For $x\in\{x_1,\dots ,x_{[{m\over 2}]}\}$, denote by $m(x,s_{\varphi })$ the multipicity of $x$ in $s_{\varphi }$ and put
$$G_x=\cases GL_{m(x;s_{\varphi })}\times GL_{m(x^{-1},s_{\varphi })}, & \hbox{\rm if $x\not\in\{\pm 1\}$},\cr GL_{m(\pm 1,s)}, & \hbox{\rm if $x=\pm 1$.} \cr\endcases$$
The group $C_{GL_N(\Bbb C)}(s_{\varphi })$ is canonically isomorphic to $\prod _xG_x$, the product going over equivalence classes of elements in the set $\{x_1,\dots ,x_{[{m\over 2}]}\}$ with respect to the relation $x\sim x^{-1}$, and to $C_{GL_N(\Bbb C)}(Im(\varphi ))$.

Denote by $H_x$ (resp. $H'_x$) the subgroup of $G_x$ defined (with $J$ an appropriate matrix which needs not to be made more precise here) by
$$\cases \{(h,Jh^{-1}J)\vert h\in GL_{m(x;s_{\varphi })}\} & \hbox{\rm if $x\not\in\{\pm 1\}$}\cr O_{m(1,s_{\varphi })}  & \hbox{\rm if $x=1$,} \cr  Sp_{m(-1,s_{\varphi })}\ \hbox{\rm (resp. $O_ {m(-1,s_{\varphi })}$),}& \hbox{\rm if $x=-1$,} \cr\endcases$$
and by $H$ (resp. $H'$) the image of $\prod _xH_x$ (resp. $\prod _xH'_x$) in $C_{GL_N(\Bbb C)}(s_{\varphi })$ by the above isomorphism.

Then, $C_{{^LG}}(Im(\varphi ))$ is isomorphic to $H$ and $C_{{^LG}}(s_{\varphi })$ to $H'$ .

In particular, $C_{C_{{^LG}}(Im(\varphi _0))}(s_{\varphi })$ and $C_{{^LG}}(Im(\varphi ))$ are only isomorphic if $m(-1,s_{\varphi })$ $=0$.

\null Proof: \rm This follows immediately from the arguments in the proof of lemma {\bf 2.2}.

\hfill{\fin 2}

\null{\bf 2.4} Remark that at the end of the following definition, we will use results from {\bf C.6} and {\bf C.7}. However, at a first reading, one may avoid to look in the annex {\bf C}.

\null{\bf Definition:} Let ${^LM}\simeq GL_{k_1}(\Bbb C)\times\cdots\times\GL_{k_d}(\Bbb C)\times {^LH_l}$ be a standard Levi subgroup of ${^LG}$ and let $\varphi _0:W_F\rightarrow{^LM}$ be a discrete normed Langlands parameter, $\gamma\mapsto (\rho _1(\gamma ),\dots ,\rho _d(\gamma ),\tau (\gamma ))$.

Recall that by {\bf 2.1}, $C_{{^LG}}(Im(\varphi _0))=\prod_{\rho } H_{\rho ;\varphi }(m(\rho ;\varphi ))$, the product going over representatives of the equivalence classes of irreducible normed representations $\rho $ of $W_F$ (w.r.t. the relation defined in {\bf 1.6}), while the $H_{\rho ;\varphi }(m)$ are complex classical groups with $H_{\rho ;\varphi }(m)$ isomorphic to $GL_m(\Bbb C)$ if $\rho $ is not self-dual, to $Sp_m(\Bbb C)$ if $\rho $ is not of type ${^LG}$, and to $O_m(\Bbb C)$ if $\rho $ is of type ${^LG}$ (with the convention $O_1(\Bbb C)=\{\pm 1\}$ if $m=1$).

Let $s$ be a semi-simple element in $C_{{^LG}}(Im(\varphi _0))$ and denote by $s_{\rho }$ the projection of $s$ on $H_{\rho ;\varphi }(m(\rho ;\varphi ))$. Define $C'_{H_{\rho ;\varphi }(m(\rho ;\varphi ))}(s_{\rho })=C_{H_{\rho ;\varphi }(m(\rho ;\varphi ))}(s_{\rho })$ except if $\rho $ and $\rho _-$ are not of the same type. In that case, put $m=m(\rho ;\varphi )$ and denote by $H'_{\rho ;\varphi }(m)$ the connected component of the $L$-group of the unramified quasi-split unitary group $U_m$ and define $C'_{H_{\rho ;\varphi }(m)}(s_{\rho })=C_{H'_{\rho ;\varphi }(m)}((-1)^{m-1}_{s_{\rho }})$ (where $(-1)_{s_{\rho }}$ is the Langlands parameter for $U_m$ defined in {\bf C.7}).

Put $C'_{C_{{^LG}}(Im(\varphi _0))}(s)=\prod _{\rho }C'_{H_{\rho ;\varphi }(m(\rho ;\varphi ))}(s_{\rho }).$ For a subset $I$ of $C'_{C_{{^LG}}(Im(\varphi _0))}(s)$, denote its centralizer by $C'_{C_{{^LG}}(Im(\varphi _0))}(s,I)$ (there is some subtlety if $\rho $ and $\rho _-$ are self-dual, but not of the same type), and write ${C'}^+_{C_{{^LG}}(Im(\varphi _0))}(s,I)$ for the subgroup of elements with determinant $1$.

\null{\bf 2.5 Lemme:} \it Let ${^LM}\simeq GL_{k_1}(\Bbb C)\times\cdots\times\GL_{k_d}(\Bbb C)\times {^LH_l}$ be a standard Levi subgroup of ${^LG}$ and let $\varphi _0:W_F\rightarrow{^LM}$ be a discrete normed Langlands parameter, $\gamma\mapsto (\rho _1(\gamma ),\dots ,\rho _d(\gamma ),\tau (\gamma ))$.

The set of equivalence classes of Langlands-Deligne parameters $\varphi :W_F\times SL_2(\Bbb C)$ $\rightarrow {^LG}$ with $\varphi _{\vert W_F}$ in the inertial orbit of $\varphi _0$ is in bijection with the set of equivalence classes of pairs $(s,\varphi _{SL_2})$ consisting of a semisimple element  $s$ and an algebraic homomorphism $SL_2(\Bbb C)\rightarrow C'_{C_{{^LG}}(Im(\varphi _0))}(s)$ by mapping $\varphi$ to $(s_{\varphi },\varphi _{\vert SL_2(\Bbb C)})$, so that $C_{\widehat{G}}(Im(\varphi ))/C_{\widehat{G}}(Im(\varphi ))^0$ is canonically isomorphic to $${C'}^+_{C_{^LG}(Im(\varphi _0))}(s_{\varphi },\varphi (SL_2(\Bbb C)))/{C'}^+_{C_{^LG}(Im(\varphi _0))}(s_{\varphi },\varphi (SL_2(\Bbb C)))^0.$$

\null Proof: \rm This is straightforward by the definitions, the above lemmas and {\bf C.7}, {\bf C.8}, remarking that $\varphi _{0,s}$ is, as element of the inertial orbit of $\varphi _0$, determined by $s_{\varphi _{0,s}}$ and that the map $s\mapsto s_{\varphi _{0,s}}$ is surjective.  \hfill{\fin 2}

\null{\bf 2.6 Theorem:} \it Let ${^LM}\simeq GL_{k_1}(\Bbb C)\times\cdots\times\GL_{k_d}(\Bbb C)\times {^LH_l}$ be a standard Levi subgroup of ${^LG}$
and let $\varphi _0:W_F\rightarrow{^LM}$ be a discrete normed Langlands parameter, $\gamma\mapsto (\rho _1(\gamma ),\dots ,\rho _d(\gamma ),\tau (\gamma ))$.

The set of equivalence classes of pairs $(s,u)$ consisting of a semisimple element  $s$ and a unipotent element $u$ in $C'_{{^LG}}(Im(\varphi _0))$ such that $sus^{-1}=u^q$ is in bijection with the set of equivalence classes of Langlands-Deligne parameters $\varphi :W_F\times SL_2(\Bbb C)\rightarrow {^LG}$ with $\varphi _{\vert W_F}$ in the inertial orbit of $\varphi _0$, so that one has a canonical isomorphism between the group of connected components of the centralizers of the images, $$C_{\widehat{G}}(Im(\varphi ))/C_{\widehat{G}}(Im(\varphi ))^0\simeq {C'}^+_{C_{^LG}(\varphi _0(W_F ))}(s,u)/{C'}^+_{C_{^LG}(\varphi _0(W_F))}(s,u)^0.$$

\null Remark: \rm Remark that $C'=C$, if $\rho _{i,-}$ is of type ${^LG}$ whenever $\rho _i$ is.

\null \it Proof: \rm  By the preceding lemma, it remains to show that the equivalence classes of pairs $(s,u)$ in the (possibly non connected) complex reductive group $C'_{{^LG}}(Im(\varphi _0))$ such that $sus^{-1}=u^q$ is in bijection with the set of equivalence classes of pairs $(s,\varphi _{SL_2})$ with $s$ in $C'_{{^LG}}(Im(\varphi _0))$ and $\varphi _{SL_2}:SL_2(\Bbb C)\rightarrow C'_{C_{{^LG}}(Im(\varphi _0))}(s)$ a morphism or algebraic groups. This is proved in \cite{KL, 2.4}. The general assumption of this paper being that the group is semi-simple and simply connected, it has been checked in \cite{H2, 3.5} that this is still valid for a connected reductive group. As Mostow's theorem is valid for possibly non-connected algebraic groups, the assumption "connected" can be relaxed, too. (The connected component of the group noted $M_{\varphi _{SL_2}}$ in  \cite{KL, 2.4} being reductive, $M_{\varphi _{SL_2}}$ is reductive.) \hfill{\fin 2}

\null{\bf 3.} In this section we will show at the end that the previous results allow to relate the category $\R _F^{\varphi _0}(G)$ to categories of unipotent representations of $p$-adic classical groups. However, before that, we will state some parameterizations of representations of collections of (possibly extended) affine Hecke algebras that follow from section {\bf 1.} and from the additional remarks about unitary groups in annex {\bf C}.

The parametrization is given by a set of conjugation classes of triples $(s,u,\Xi )$ associated to a given complex group, where $s$ is a semi-simple element, $u$ a unipotent element such that $sus^{-1}=u^q$ and $\Xi $ an irreducible representation of the group of components of the common centralizer of $s$ and $u$.

\null{\bf 3.1} The following is the special case of {\bf 1.7}, {\bf 1.8} for $G=SO_{2d+1}$, $M=T$ and $\rho $ the trivial representation that is treated in \cite{L2} with modifications in \cite{W}.

\null{\bf Theorem: } \it Fix an integer $d\geq 1$. If $(d_+, d_-)\ne 0$ is a pair of integers which are each one products of two consecutive integers, $d_+=a_+(a_++1)$ and $d_-=a_-(a_-+1)$, with $d_++d_-\leq 2d+1$, denote by $\H (d_+, d_-)$ the affine Hecke algebra with root datum equal to the one of $SO_{2d+1-d_+-d_-}$ and unequal parameters $q,\dots ,q, q^{a_++a_-+1}; q^{\vert a_+-a_-\vert}$. Denote by $\H (0,0)$ the affine Hecke algebra with root datum equal to the one of $Sp_{2d}(\Bbb C)$ and equal parameters $q,\dots ,q$.

Then, the set of triples $(s,u,\Xi )$ associated to the group $Sp_{2d}(\Bbb C)$ with $\Xi (-1)=1$ is in natural bijection with the set $$\bigcup_{(d_+,d_-), {d_++d_-\over 2}\ even}(irreducible\ right-\H(d_+, d_-)-modules)$$
and the one with $\Xi (-1)=-1$ is in nature bijection with the set $$\bigcup_{(d_+,d_-), {d_++d_-\over 2}\ odd}(irreducible\ right-\H(d_+, d_-)-modules).$$ \rm

\null\it Remark: \rm By natural bijection, we mean what is implied by the properties of the local Langlands correspondence stated in {\bf 1.7}. We will not explain this here more, except that compact $s$ should correspond to tempered representations, discrete $(s,u)$ (i.e. those which do not lie in a proper parabolic subgroup) to discrete series representations and that Langlands-Shahidi local factors (defined for Hecke-algebra representations by equivalence of categories) equal the Artinian ones (deduced from $(s,u)$). The bijection associates implicitly to a triple $(s,u,\Xi )$ a "supercuspidal support" according to {\bf 1.6}, but we will not give a precise construction here, neither in the subsequent cases. (For the above case, it can be found in \cite{L2} and \cite{W}.)

\null{\bf 3.2} The following follows by combining {\bf 1.7} and {\bf 1.8} to the special case $G=Sp_{2d}(F)$ (for (i)) and $G=O_{2d}(F)$ (for (ii)) with $\rho $ the trivial representation and $M$ the maximal split torus. Remark that these groups are not of adjoint type, so that \cite{L2} does not apply to these cases, but a slight generalization should do. However, this has not been written yet.

\null {\bf Theorem:} \it Fix an integer $d\geq 1$.

If  $(d_+, d_-)\ne 0$ is a pair of integers, which are squares, such that $d_++d_-\leq 2d+1$, denote by $\H (d_+, d_-)$ the affine Hecke algebra with root datum equal to the one of $SO_{2d+1-d_+-d_-}$, if  $d_++d_-$ is even, and equal to the one of $SO_{2d+2-d_+-d_-}$, if $d_++d_-$ is odd,and unequal parameters $q,\dots ,q, q^{\sqrt{d_+}+\sqrt{d_-}}; q^{\vert \sqrt{d_+}-\sqrt{d_-}\vert}$.

In addition, denote by $\H (0,0)$ the semi-direct product of an affine Hecke algebra with equal parameter $q$ and root datum equal to the one of $SO_{2d}$ with the group algebra of a finite cyclic group of order 2, which acts by the outer automorphism of the root system. Define
$\epsilon^+ _{d_+,d_-}=\cases 4 &\hbox{\rm if $d_+$ even, $d_++d_-\in 8\Bbb Z$, $d_+\cdot d_-\ne 0$,}\cr 1 &\hbox{\rm if $d_+=d_-=0$,}\cr 0 &\hbox{\rm if $d_+$ even, $d_++d_-\in 4\Bbb Z\setminus8\Bbb Z$,} \cr 2&\hbox{\rm otherwise.}\cr\endcases$ and put $\epsilon^- _{d_+,d_-}=4-\epsilon^+ _{d_+,d_-}$ if $d_+\cdot d_-\ne 0$, $\epsilon^- _{d_+,d_-}=2-\epsilon^+ _{d_+,d_-}$ if exactly one of $d_+$ and $d_-$ is $0$, and $\epsilon^0_{0,0}=0$.

\null i) Denote by $S_o$ the set of pairs of integers $(d_+, d_-)$ such that $d_+$ and $d_-$ are squares, $d_++d_-$ is odd and $\leq 2d+1$. (Consequently $\epsilon ^\pm _{d_+,d_-}=2$.
The set of triples $(s,u,\Xi )$ associated to the group $SO_{2d+1}(\Bbb C)$ is in natural bijection with the multiset $$\bigcup_{(d_+,d_-)\in S_o} 2\ (irreducible\ right-\H(d_+, d_-)-modules).$$

\null ii) Denote by $S_e$ the set of pairs of integers $(d_+, d_-)$ such that $d_+$ and $d_-$ are squares, $d_++d_-$ is even and $\leq 2d+1$.
The set of triples $(s,u,\Xi )$ associated to the group $O_{2d}(\Bbb C)$ with $\Xi _{\vert Z_{\widehat{G}}}=\pm 1$ is in natural bijection with the multiset $$\bigcup_{(d_+,d_-)\in S_e} \epsilon^\pm _{d_+,d_-}\ (irreducible\ right-\H(d_+, d_-)-modules).$$ \rm

\null{\bf 3.3} Assuming the appropriate parts of the Langlands correspondence {\bf 1.7} established for the non split pure inner form of the even unramified quasi-split unitary group, the following follows by combining {\bf 1.7} and {\bf 1.8} as generalized to quasi-split unitary groups in annex {\bf C.}. Here $U_m$  will denote the unramified quasi-split unitary group of semi-simple rank $m$.

Fix $m$. A triple $(s,u,\Xi )$ as above relative to $^LU_m$, will be said associated to $^LU_m$, if $s$ is not an element of the connected component $({^LU_m})^0$ of $^LU_m$.

\null{\bf Theorem:} \it If $d_+$ is a square integer and $d_-$ the product of two consecutive integers, $d_-=a_-(a_-+1)$, denote by $\H(d_+,d_-)$ the affine Hecke algebra with root datum equal to the one of $SO_{2d+1-d_+-d_-}$, if $d_++d_-$ is even, and to the one of $SO_{2d+2-d_+-d_-}$, if $d_++d_-$ is odd, and unequal parameters $q,\dots ,q, q^{\sqrt{d_+}+a_-+{1\over 2}};$ $q^{\vert \sqrt{d_+}-a_--{1\over 2}\vert}$.

In addition, denote by $\H (0,0)$ the affine Hecke algebra with root datum equal to the one of $SO_{2d+1}$ and unequal parameters $q,q,\dots ,q,q^{1/2}$.

Denote by $S_e$ the set of pairs $(d_+, d_-)$ with $d_+$ an even square, $d_-$ the product of two consecutive integers, $d_++d_-\leq 2d+1$, and by $S_o$ the set of pairs $(d_+, d_-)$ with $d_+$ an odd square, $d_-$ the product of two consecutive integers, $d_++d_-\leq 2d+1$.

(i) If $m$ is an odd integer, $m=2d+1$, then the set of triples $(s,u,\Xi)$ associated to $^LU_m$ with $\Xi_{\vert \{\pm 1\}}$ fixed is in natural bijection with the multiset
$$\bigcup_{(d_+,d_-)\in S_o} (irreducible\ right-\H(d_+, d_-)-modules).$$

(ii) If $m$ is an even integer, $m=2d$, put $\epsilon _{d_+,d_-}=2$, if $d_+\ne 0$, and $\epsilon _{d_+,d_-}=1$ otherwise. Then the set of triples $(s,u,\Xi)$ associated to $^LU_m$ with $\Xi(-1)=1$ is in natural bijection with the multiset
$$\bigcup_{(d_+,d_-)\in S_e, {d_+\over 2}\ even} \epsilon _{d_+,d_-}\ (irreducible\ right-\H(d_+, d_-)-modules)$$ and the one with $\Xi(-1)=-1$ is in natural bijection with the multiset $$\bigcup_{(d_+,d_-)\in S_e, {d_+\over 2}\ odd} 2\ (irreducible\ right-\H(d_+, d_-)-modules). $$\rm

\null{\bf 3.4} If $t_{\rho }$ is an integer $\geq 1$, denote by $F_{t_{\rho }}$ the unramified extension of $F$ of degree $t_{\rho }$, which is unique in a given algebraic closure of $F$. If $\varphi $ is a Langlands parameter which is not normed and $\varphi _0$ is the normed Langlands parameter in its orbit, we put $\R^{\varphi }=\R^{\varphi _0}$. We also put $\varphi _0=1$, if $\varphi _0$ is the Langlands parameter relative to the minimal standard Levi subgroup that corresponds to the trivial representation.

\null {\bf Theorem:} \it Assume that there is an irreducible representation $\rho $ such that all irreducible components of $\iota\circ\varphi _0$ are either isomorphic to $\rho $ or to $\rho ^{\vee }$. Then, with $m=m(\rho ;\varphi )$, one has

(i) if $\rho $ is not self-dual, then the category $\R_F^{\varphi _0}(\underline{G})$ is equivalent to $\R _{F_{t_{\rho }}}^{1}(GL_m)$.

(ii) if $\rho $ is self-dual and not of type $^LG$, then the category $\R_F^{\varphi _0}(\underline{G})$ is equivalent to $\R _{F_{t_{\rho }}}^1(SO_{m+1})$.

(iii) if $\rho $ and $\rho _-$ are both of type $^LG$, then the category $\R_F^{\varphi _0}(\underline{G})$ is equivalent to $\R _{F_{t_{\rho }}}^1(Sp_{m-1})$ if $m$ is odd, and to $\R _{F_{t_{\rho }}}^1(O_m)$ otherwise.

(iv) if $\rho $ and $\rho _-$ are self-dual but not of the same type, then, with $U_m$ equal to the unramified quasi-split unitary group of absolute rank $m$, the category $\R_F^{\varphi _0}(\underline{G})$ is equivalent to $\R _{F_{t_{\rho }}}^1(U_m)$.

The same holds, if one replaces $\R ^\cdot$ by $\R^{\cdot,+}$ or $\R^{\cdot,-}$.

\null Proof: \rm This follows from theorem {\bf 1.9} (and its generalization to unitary groups in {\bf C.5} together with proposition {\bf C.6}) applied to the above cases for $\varphi _0$, after remarking that in each of the cases the sets $\S^\cdot (\varphi _0)$ and $\S^\cdot (1)$ are equal, while the alternating characters associated to their elements are the same. There is no need here to restrict to finitely generated representations, as the tensor product is not involved. \hfill{\fin 2}

\null{\bf 3.5} Recall the equivalence relation on $supp(\varphi _0)$ given by $\rho\sim\rho^{\vee }$ introduced in {\bf 1.6} and that the index $f$ denotes the full subcategory of finitely generated representations.

\null{\bf Corollary:} \it  The category $\R^{\varphi _0}_F(\underline{G})_f$ is equivalent to $$\bigotimes _{\rho\in supp(\varphi _0)/\sim }\R _{F_{t_{\rho }}}^1(H_{\rho }(m(\rho ;\varphi _0)))_f$$ with $H_{\rho }(m)$ equal to $GL_m$, $SO_{m+1}$, $Sp_{m-1}$, $O_m$ or the unramified quasi-split unitary group $U_m$, if respectively $\rho $ is not self-dual, not of type $^L{G}$, $\rho $ and $\rho _-$ are both of type $^L{G}$ with $m$ odd, with $m$ even, or $\rho $ and $\rho _-$ are self-dual but not of the same type.

The same holds, if one replaces $\R ^\cdot$ by $\R^{\cdot,+}$ or $\R^{\cdot,-}$.

\null\it Remark: \rm As follows from remark (iii) after theorem {\bf 1.8}, this equivalence of category preserves temperedness. Discreteness is preserved if none of the $H_{\rho }$ is a general linear group - otherwise there are no discrete series representations in $\R^{\varphi _0}_F(\underline{G})_f$. Unitarity is conjectured.

\null \it Proof:  \rm Denote by  $\S (\varphi _0)_{\rho }$ the projection of $\S (\varphi _0)$ on the $\rho $'s component and by $\varphi _{0,\rho }$ the discrete Langlands parameter (unique up to equivalence) into some $L$-group $^LH_{\rho }$ of the same type as $^LG$ that satisfies the following condition (which determines also the semi-simple rank of $^LH_{\rho }$)
$$\iota\circ\varphi _{0,\rho }=\cases m(\rho;\varphi _0)\rho, & \hbox{\rm if $\rho\simeq\rho^{\vee }$};\cr m(\rho;\varphi _0)\rho\oplus m(\rho^{\vee };\varphi _0)\rho^{\vee }, & \hbox{\rm otherwise.} \cr\endcases$$

Then one has $\S(\varphi _0)=\prod _{\rho }\S(\varphi _0)_{\rho }$, $\S(\varphi _0)_{\rho }=\S(\varphi _{0,\rho })$. In addition, if $S=(S_{\rho })_{\rho }\in\S(\varphi _0)$, then $\H_{\varphi _0,S,\rho  }=\H_{\varphi _{0,\rho },S_{\rho },\rho  }$ and the group of alternating characters $\widehat{S}$ is the product of the groups $\widehat{S_{\rho }}$, where $S_{\rho }$ is seen as element of $\S(\varphi _{0,\rho })$.

Applying {\bf B.4} to {\bf 1.9}, one gets from this
$$\eqalign{\R_F^{\varphi _0}(\underline{G})_f\simeq&\bigoplus_{S\in \sS(\varphi _0), \epsilon\in \widehat{S}^\pm} (\bigotimes _{\rho\in supp(\varphi _0)/\sim} (right-\H_{\varphi _0,S,\rho  }-modules)_f\ )\cr \simeq &\bigotimes _{\rho\in supp(\varphi _0)/\sim} (\bigoplus_{S_{\rho }\in \sS(\varphi _{0,\rho }), \epsilon\in \widehat{S}_{\rho }^\pm} (right-\H_{\varphi _{0,\rho },S_{\rho },\rho  }-modules)_f).\cr}$$
Using {\bf 3.4} and {\bf B.5}, the statement of the corollary follows. \hfill{\fin 2}

\null  \it Remark: \rm A statement of this kind had been conjectured by G. Lusztig \cite{L4, section 19}. Once the results in \cite{L2} appropriately generalized to symplectic, unitary and the (full) even orthogonal group, one should be able to describe multiplicities in standard modules from intersection cohomology as described in \cite{L4, section 19}.

\null
\null

{\bf A. Annex:} Equivalence of categories for the full orthogonal group

\null {\bf A.1} The aim of this annex is to show how the results of \cite{H3, H4} generalize to the full orthogonal group, which is not connected. So, in this annex, $H$ will denote a pure inner form of a full split orthogonal group, either split or not. The case when its connected component $H^0$ is an odd orthogonal group is quite easy. Then $H$ is isomorphic to a direct product $H^0\times\{\pm 1\}$. The Levi subgroups of $H$ are of the form $M=M^0\times\{\pm 1\}$, where $M^0$ is a Levi subgroup of $H^0$, so that the supercuspidal representations of $M$ are of the form $\sigma^0 \eta $, where $\sigma ^0$ is a supercuspidal representation of $M^0$ and $\eta $ a character of $\{\pm 1\}$. One sees immediately that the restriction to $\{\pm 1\}$ of a representation in the supercuspidal support of an irreducible representation $\pi $ of $H$ is determined by the restriction of $\pi $ to this group. So, one may decompose $Rep(H)$ as a direct sum of subcategories $Rep_{M^0,\so }(H^0)\oplus Rep_{M^0,\so ,-1}(H^0)$, where the $Rep_{M^0,\so }(H^0)$ denote the Bernstein components for $H^0$ and the $Rep_{M^0,\so ,-1} (H^0)$ denote the part with non-trivial restriction to $\{\pm 1\}$. As the results of \cite{H3,H4} apply to $Rep(H^0)$, we are done.

\null{\bf A.2} Assume now for the rest of this annex that $n$ is even. Then, $H$ is isomorphic to a semi-direct product $H^0\rtimes\{1,r_0\}$, where $H^0$ is an even orthogonal group and $r_0$ is of order $2$ and acts on $H^0$ by the outer isomorphism. We refer to \cite{GH} for results for the representation theory of a non connected reductive group. We consider only Levi subgroups which are \it cuspidal \rm in the terminology of \cite{GH}. In particular, one deduces from this paper that the Bernstein decomposition is still valid and that, if $M$ is a Levi subgroup of $H$ and $\o $ denotes the inertial orbit of an irreducible supercuspidal representation of $M$, then, with the notations in \cite{H3}, $i_P^HE_{B_{\so }}$ is a projective generator of $Rep_{(M,\so )}(H)$, which implies that the category $Rep_{(M,\so )}(H)$ is equivalent to the category of right-modules over $End_H(i_P^HE_{B_{\so }})$ by Morita theory.

The aim of this annex is to show that $End_H(i_P^HE_{B_{\sso }})$ has the form given in theorem {\bf 1.8}.

\null {\bf A.3} Denote by $W^0$ the Weyl group of $H^0$ and define $W:=W^0\rtimes\{1,r_0\}$, and similar for the Weyl group $W^M$  of a Levi subgroup $M$ of $H$. If $M$ is a Levi subgroup of $H$, define $W(M)=W^M\backslash\{w \in W\vert w^{-1}Mw=M\}$ and similarly $W^0(M^0)$, which will also be denoted (abusively) $W^0(M)$.

\null {\bf Lemma:} \it One has $W(M)=W^0(M)$, except if $M$ is isomorphic to a product of general linear groups and at least one of them has odd rank. In particular, $W(M)=W^0(M)$ if $H$ is not quasi-split.

\null Proof: \rm If $M$ has a factor $H_l$ with $l\geq 2$, then $r_0\in M$. So, $M$ has to be a product of linear groups if $W(M)\ne W^0(M)$. If $M$ is a product of general linear groups of even rank, then every element of $w$ that satisfies $w^{-1}Mw=M$ must have an even number of sign changes $x\mapsto x^{-1}$ on the maximal torus. This means that it lies in $W^0$. If $M$ is a product of general linear groups, one of them being of odd rank $k$, one sees that there is an element in $W$ which induces the outer automorphism on $GL_k$ and which does not lie in $W^0$. \hfill{\fin 2}

\null {\bf A.4} Let $\o $ be the inertial orbit of a supercuspidal representation of a Levi subgroup $M$ of $H$. Its restriction to $M^0$ decomposes into one or two inertial orbits. Fix an orbit $\o ^0$ in the restriction. Denote by $W(M,\o )$ (resp. $W^0(M,\o )$) the subset of elements of $W(M)$ (resp. $W^0(M)$) which stabilize $\o $ (resp. $\o ^0$).

\null {\bf Lemma:} \it One has $W(M,\o )=W^0(M,\o )$ except if $M$ is a product of general linear groups and at least one factor of $\o ^0$ is the inertial orbit of a self-dual representation of a general linear group of odd rank.

\null Proof: \rm The group $W(M,\o )$ is a subgroup of $W(M)$. It follows that the equality $W(M,\o )=W^0(M,\o )$ can only fail if $M$ is a product of general linear groups and at least one of them has odd rank $k$. In addition, one sees that at least one factor of $\o $ corresponding to a $GL_k(F)$ with $k$ odd must be the orbit of a self-dual representation. \hfill{\fin 2}

\null {\bf A.5} Denote by $R(\o )$ the subgroup of elements $r$ of $W(M,\o )$ that send positive roots for $M$ to positive roots. Define $R^0(\o )=R(\o )\cap W^0(M,\o )$. Recall \cite{H3} that $W^0(M,\o )$ is a semi-direct product $W_{\so }^0\rtimes R^0(\o )$, so that one has $W(M,\o )=W^0_{\so }\rtimes R(\o )$. As $ind^M_{M^0}E_{B_{\sso ^0}}$ is either equal to $E_{B_{\sso }}$ or a direct sum $E_{B_{\sso }}\oplus E_{B_{\sso' }}$, one can define, for $w\in W^0(\o ^0)$ and $r\in R^0(\o )$, operators $T_w$ and $J_r$ in $End_H(i_P^HE_{B_{\sso }})$ from the ones for $End_{G^0}(i_{P^0}^{H^0}E_{B_{\sso ^0}})$ by induction. If $r\in R(\o )\setminus R^0(\o )$, note $\lambda (r)$ the action of $r$ on $i_P^HE_{B_{\sso }}$ by left-translation and by $\tau _r$ the one of $r$ on $B_{\so }$ by right translation \cite{H3}, and put $J_r=\tau _r\lambda (r)$. These operators $J_r$ commute obviously with the other operators $J_{r'}$, $r'\in R(\o )$, and satisfy the commuting relation $T_wJ_r=J_rT_{r^{-1}wr}$ for $w\in W_{\so }^0$.

\null {\bf Lemma:} \it The operators $sp_{\chi }J_rT_w$, $r\in R(\o )$, $w\in W(M,\o )$, are linearly independent for all $\chi\in\X(M)$.

\null Proof: \rm The proof of \cite{H3, 5.9} generalizes, as the commuting relations for the operators $J_r$, $r\in R(\o )$ are still the same. \hfill{\fin 2}

\null {\bf A.6 Lemma:} \it One has $\Hom _H(i_P^HE_{B_{\sso }}, i_P^GE_{K(B_{\sso })})=\bigoplus _{w,r}K(B_{\so })J_rT_w.$

\null Proof: \rm This follows from the linear independence and the computation of the Jacquet module with help of the geometric lemma in the non connected case \cite{C, 4.1}, taking into account lemma {\bf A.4}. \hfill{\fin 2}

\null{\bf A.7 Theorem:} \it  One has $\End _H(i_P^HE_{B_{\sso }})=\bigoplus _{w,r}B_{\so }J_rT_w.$

\null Proof: \rm The proof of \cite{H3, 5.10} generalizes, as the commuting relations for the operators $J_r$, $r\in R(\o )$ are still the same. \hfill{\fin 2}

\null \it Remark: \rm  As the $T_w$ satisfy the same relations as their restrictions to the space of the representation of the connected component, it follows that $End_H(i_H^GE_{B_{\sso }})$ is an (possibly extended) affine Hecke algebra isomorphic to $End_{H^0}(i_{P^0}^{G^0}E_{B_{\sso ^0}})$, except if $M$ is a product of general linear groups and at least one factor of $\o $ is the inertial orbit of a self-dual representation of a general linear group of odd rank. In this case, one has additional operators $J_r$ with $r\in R(\o )\setminus R^0(\o )$ as stated in {\bf 1.8}.

In fact, we have an erratum to \cite{H3, H4} w.r.t. the statements for the even dimensional special orthogonal group, the connected component of $H$: if $M=M^0$ (i.e. $M$ is a product of general linear groups), $End_{H^0}(i_{P^0}^{H^0}E_{B_{\sso ^0}})$ is in general isomorphic to a tensor product $\bigotimes _{\rho }\H_{\rho }\otimes((\otimes _{\rho '}\H_{\rho '}^0)\rtimes\Bbb C[R_{nq}]),$ the first product going over elements $\rho $ in the support of the normed Langlands parameter $\varphi_0$ associated to $\o $ which are not odd orthogonal and the second product over the odd orthogonal ones, $R_{nq}$ being generated by Weyl group elements that send positive roots in $\Sigma _{\so }$ (in the notations of \cite{H3}) to positive roots and have sign changes $x\mapsto x^{-1}$ on two factors $GL_k(F)$ and $GL_{k'}(F)$, on which odd orthogonal representations with distinct inertial orbits are defined. Here the $\H_{\rho }$ denote the (possibly extended) affine Hecke algebras from {\bf 1.8} and $\H_{\rho }^0$ the affine Hecke algebra part (i.e. omitting the finite group algebra part, if there is any). One remarks that the above semi-direct product is with a tensor product of affine Hecke algebras associated to odd orthogonal representations in the support, but \it does not \rm decompose into a tensor product of semi-direct products of the different affine Hecke algebras with a group algebra.

\null\null{\bf Annex B:} Tensor product of abelian categories

\null {\bf B.1} \it Definition \rm \cite{D, 5.} Let $k$ be a commutative ring and $(\hbox{\main A }_i)_{i\in I}$  a finite family of $k$-linear abelian categories. A $k$-linear abelian category $\hbox{\main A}$ equipped with a $k$-multilinear functor right exact in each variable $$\otimes: \prod \A _i\rightarrow \A$$
is called \it tensor product over $k$ \rm of the categories $\hbox{\main A }_i$ if and only if the following condition is satisfied: denote for a $k$-linear abelian category $\hbox{\main C}$ by $\underline{Hom}_{\hbox{\sit k,e\ \`a\ d}}(\A, \C)$ the category of right exact functors from $\A $ to $\C $ and by $\underline{Hom}_{\hbox{\sit k,e\ \`a\ d}}((\A )_{i\in I}, \C)$ the category of right exact functors multilinear in each variable from the product of the $\A _i$ to $\C $.

One asks then that for every category $\C $ the composed functor with the above
$$\underline{Hom}_{\hbox{\sit k,e\ \`a\ d}}(\A, \C)\rightarrow \underline{Hom}_{\hbox{\sit k,e\ \`a\ d}}((\A _i)_{i\in I}, \C)$$ is an equivalence of categories.

\null {\bf B.2 Proposition:} \cite{D, 5.3} \it Let $(A _i)_{i\in I}$ be a finite family of coherent $k$-algebras that have a coherent tensor product over $k$. Denote by $(A _i)_{coh}$ (resp. $(\otimes A_i)_{coh}$) the corresponding abelian category of right modules of finite presentation. The tensor product over $k$
$$\otimes:\prod (A_i)_{coh}\rightarrow (\otimes A_i)_{coh}$$ defines $(\otimes A_i)_{coh}$ as tensor product over $k$ of the $(A_i)_{coh}$. \rm

\null {\bf B.3 Proposition:} \it (i) An extended affine Hecke algebra with unequal parameters is a coherent $\Bbb C$-algebra.

(ii) Any finitely generated right module over an extended affine Hecke algebra with unequal parameters is coherent.

\null Proof: \rm (i) An affine Hecke algebra with unequal parameters is a free module of finite rank over the group ring of a finitely generated lattice. As the group ring of a finitely generated lattice is noetherian as quotient of a polynomial ring, the extended affine Hecke algebra is noetherian as a module. But every ideal of this algebra is a submodule. So, it is finitely generated. One concludes that an affine Hecke algebra is noetherian and in particular coherent. As an extended affine Hecke algebra is, as a module, the sum of an affine Hecke algebra with a finite dimensional $\Bbb C$-vector space, we are done.

(ii) A finitely generated right-module over a noetherian $\Bbb C$-algebra is coherent.

\hfill{\fin 2}

\null{\bf B.4 Proposition:} \it Let $k$ be a commutative ring and $(\hbox{\main A }_i)_{i\in I}$  a finite family of $k$-linear abelian categories. Assume that each $\A _i$ is a direct sum of $k$-linear categories $\A _{i,j}$, $j=1,\dots ,l_i$. Suppose that for each family of integers $\underline{j}=(j_i)_{i\in I}$, $1\leq j_i\leq l_i$, the family $(\A _{i,j_i})_{i\in I}$ has a tensor product $\A_{\underline{j} }$. Then, the tensor product of the categories $\A _i$ is isomorphic to the  direct sum of the categories $\A_{\underline{j}}$.

\null Proof:  \rm Write $\J$ for the set of the $\underline{j}$. One has an equivalence of categories between $\prod_i (\bigoplus _{j=1}^{l_i}\A _{j,i})$ and $\bigoplus_{\underline{j}\in\sJ}\prod _i\A _{j_i,i}$, and consequently between $\underline{Hom}_{\hbox{\sit k,e\ \`a\ d}}$

\noindent{$((\bigoplus _{j=1}^{l_i}\A _{j,i})_{i\in I}, \C)\ \ $   and $\bigoplus_{\underline{j}\in\sJ}\underline{Hom}_{\hbox{\sit k,e\ \`a\ d}}((\A _{j_i,i})_{i\in I}, \C)$. Denote by $\A _{\underline{j}}$ the tensor product of  $(\A _{j_i,i})_{i\in I}$. One sees immediately that $\bigoplus_{\underline{j}\in\sJ} \A _{\underline{j}}$ satisfies the universal property for the tensor product of $(\bigoplus _{j=1}^{l_i}\A _{j,i})_{i\in I}$.} \hfill{\fin 2}

\null {\bf B.5 Proposition:} \it Let $(\H _i)_{i\in I}$ be a finite family of extended affine Hecke algebras with parameters. Let $\B _i$ be a finite family of $k$-linear abelian categories with each $\B _i$ equivalent to the category $(\H_i)_f$ of finitely generated modules over $\H _i$. Then, the tensor product of the $k$-linear abelian categories $\B_i $ exists and is equivalent to the tensor product of the categories $(\H_i)_f$.

\null Proof: \rm The equivalence of categories $\B _i\rightarrow (\H_i)_f$ gives equivalences of categories $\prod \B _i\rightarrow \prod (\H_i)_f$ and $\underline{Hom}_{\hbox{\sit k,e\ \`a\ d}}(((\H_i)_f)_{i\in I}, \C)\rightarrow \underline{Hom}_{\hbox{\sit k,e\ \`a\ d}}((\B_i )_{i\in I}, \C)$. With this, it is immediate that the $(\B_i)_{i\in I}$ satisfy the universal property with respect to the tensor product of the $\H _i$. \hfill{\fin 2}

\null\null{\bf Annex C:} The case of the unitary group

\null\null {\bf C.1} In this annex, we will show that the results of the section {\bf 1.} and {\bf 2.} generalize to quasi-split unitary groups. We will give a few remarks, justifying that \cite{H3} generalizes to pure inner forms of unitary groups. To obtain the generalization of  \cite{H4}, the reference to \cite{M1} has to be replaced by \cite{M2} (see also \cite{M3}). (Remark that Arthur's work has been generalized to quasi-split unitary groups in \cite{Mk}. Based on that, inner forms of unitary groups are treated in \cite{TMSW}.)

To generalize section {\bf 2.}, we rely on \cite{GGP} for appropriate results for Langlands parameters for unitary groups.

\null {\bf C.2} In this section, $H$ will denote the group of $F$-rational points of a quasi-split unitary group $\underline{H}$ with respect to a quadratic extension $E/F$. As $\ul{H}$ is not split,  the $L$-group of $H$ is a semi-direct product $GL_n(\Bbb C)\rtimes Gal(E/F)$, where $GL_n(\Bbb C)$ is the Langlands dual group of $H$.

According to the parity of $n$, we will say that $H$ is an even or odd unitary group. We will denote by $W_E$ the Weil group of $E$. The notion of a conjugate-orthogonal and a conjugate-symplectic representation of $W_E$ is defined in \cite{GGP}. A conjugate-dual representation $\rho $ of $W_E$ will be said of type $^LH$ if either $n$ is even and $\rho $ is conjugate-sympletic, or $n$ is odd and $\rho $ is conjugate-orthogonal. Otherwise, we will say that $\rho $ is not of type $^LH$. We stretch that the use of either these notions will presume that $\rho $ is conjugate-dual. The same terminology will also be used when $W_E$ is replaced by the Weil-Deligne group $W_E\times SL_2(\Bbb C)$.

There is a unique pure inner form of $H$ which we will denote by $H^-$. If $n$ is odd, then $H^-$ is isomorphic to $H$ and if $n$ is even then $H^-$ is not quasi-split. We will write again sometimes $H^+$ for $H$.

\null {\bf C.3} A Langlands parameter for $H$ is a morphism $W_F\rightarrow\ ^LH$ such that the projection to the first factor is a Langlands parameter (as defined in {\bf 1.}) and the projection to the second factor is the projection $W_F\rightarrow Gal(E/F)$. The definition of a Langlands-Deligne parameter is straighforward.

It is explained in \cite{GGP, section 8} that Langlands and Langlands-Deligne parameters for $H$ are in bijective correspondence with conjugate-dual representations of type $^LH$ of $W_E$ or $W_E\times SL_2(\Bbb C)$ respectively. It follows from this also that it does not matter to define equivalence for Langlands parameters or Langlands-Deligne parameters by conjugation by an element of $\widehat{H}$ or $^LH$. If $\varphi $ is a Langlands or a Langlands-Deligne parameter for $H$ we will denote by $\varphi _E$ the corresponding conjugate-dual representation of type $^LH$.

With this terminology, replacing $\iota \circ\varphi $ by $\varphi _E$, it is shown in \cite{M2, 8.4.4} (see also \cite{M3}) that the part of theorem {\bf 1.1} that applies to $H^+$ generalizes. As $H$ is isomorphic to $H^-$ in the odd case, one sees easily that this implies the whole theorem {\bf 1.1} in the odd case. As \cite{M2, M3} do not treat the pure inner form of the even quasi-split unitary group, this case has to be left as a conjecture, but it is certainly true.

\null{\bf C.4} The definition in {\bf 1.2} has to be modified to choose in each inertial class of an irreducible representation of $W_E$ a base point that is conjugate-dual if there is such a representation in the inertial class and, if possible, even conjugate dual of the same type than $^LH$.

A standard Levi subgroup $M$ of $H$ has the form $GL_{k_1}(E)\times\cdots\times GL_{k_r}(E)\times H_l$, where $H_l$ is a unitary group of the same type (even or odd) than $H$. One has the equality $n=2(k_1+\cdots +k_r)+L$, where $L$ is defined by $\widehat{H_l}=GL_L(\Bbb C)$.  One has then $^LM=GL_{k_1}(\Bbb C)\times\cdots\times GL_{k_r}(\Bbb C)\times\ ^LH_l$. If $\varphi =(\rho _1,\dots ,\rho _k,\tau ) : W_F\rightarrow\ ^LM$ is a discrete Langlands parameter, we will denote by $\rho _{i,E}$ the corresponding irreducible representation $W_E\rightarrow GL_{k_i}(\Bbb C)$, by $^c\rho _{i,E}$ the conjugate representation and by $\varphi _E$ the conjugate-dual representation $W_E\rightarrow  GL_N(\Bbb C)$ of type $^LH$  that is isomorphic to $\tau _E\oplus\bigoplus_{i=1}^k(\rho _{i,E}\oplus\  ^c\rho _{i,E}^{\vee })$. We  will call $\varphi $ normed, if $\varphi _E$ is normed in the sense defined by {\bf 1.2}.

If $s$ is an element in the centralizer of $\varphi _E(W_F)$ in $\widehat{G}$ such that the representation $\varphi _{E,s}$ in the inertial class of $\varphi _E$ is conjugate-dual of type $^LG$, then we will denote by $\varphi _s$ the corresponding Langlands parameter for $^LG$. The set of the $\varphi _s$ will be the \it inertial orbit \rm of $\varphi $. The proof of proposition {\bf 1.3} generalizes, after replacing $\iota\circ\varphi $ by $\varphi _E$ and remarking that $W_E\cap I_F=I_E$, and one sees that two Langlands parameters relative to $H$ lie in the same inertial orbit, if and only if their restriction to the inertia subgroup of $I_F$ are conjugated by an element of $^LH$. One defines the multiplicity $m(\rho ;\varphi )$ to be the multiplicity of $\rho $ in $\varphi _E$.

The proposition {\bf 1.4} generalizes obviously to representations of $W_E$, replacing self-dual by dual-conjugate, remarking that $\vert\cdot\vert _E$ is self-conjugate. One defines then for a conjugate-dual representation $\rho $ the representation $\rho _-$ accordingly.

Replacing self-dual by dual-conjugate, the generalizations of the notions defined in {\bf 1.6} is straightforward and the theorem at the end remains valid.

The definition of the category $Rep_F^{\varphi _0}(\underline{H})$, for $\varphi _0$ a normed Langlands parameter for $^LH$, and its subcategories $Rep_F^{\varphi _0,\pm }(\underline{H})$ in the statement of the local Langlands correspondence {\bf 1.7} is then clear. (The cases with index "-" do not appear, i.e. these notions can remain undefined.) The $L$-functions and local factors which have to be used here are those coming from the Asai representation.

\null{\bf C.5} The theorem {\bf 1.8} is based on \cite{H3} and \cite{H4}, which do not explicitly include unitary groups. However, \cite{H3} generalizes with only minor changes to pure inner forms of quasi-split unitary groups: as the Levi subgroups of $H$ are of the form $GL_{k_1}(E)\times\cdots\times GL_{k_r}(E)\times H_l$, the assumptions made in \cite{H3, 1.13 - 17} and the results therein remain valid when taken the absolute value and a uniformizer for $E$ at appropriate places. One remarks that the relative reduced roots for $H$ form a root system of type $B$ in the odd case and of type $C$ in the even case. From this, the generalization of \cite{H3, 1.13} is straightforward. The same applies to section {\bf 6.} and {\bf 7.} of \cite{H3}.

The Plancherel measure of a representation of type $\sigma _{S,\epsilon }$ of a Levi subgroup $M_S$ can be computed as in \cite{H4} according to the results in \cite{M2, M3} (especially \cite{M2, 8.4.4} already mentioned above in {\bf C.3}), the relation with reducibility points remaining the same as in the orthogonal or symplectic case, using $\vert\cdot\vert _E$
instead of $\vert\cdot\vert _F$. Replacing self-dual by conjugate-dual, the generalization of {\bf 1.8} is straightforward. The corollary {\bf 1.9} is then a direct consequence.

\null{\bf C.6 Proposition:} \cite{GGP, 3.4} \it The trivial character of $E^{\times }$ is always a conjugate-orthogonal representation. The nontrivial unramified quadratic character of $E^{\times }$ is conjugate-symplectic, if and only if $E/F$ is unramified. Otherwise, it is conjugate-orthogonal. \rm

\null {\bf C.7} The unitary group $H$ is called unramified if $E/F$ is an unramified extension. Denote by $1$ the Langlands parameter for $H$ such that $1_E$ is $n$ times the trivial representation of $W_E$. We will write $-1$ for the Langlands parameter $1_{-1}$ for $H$ in the above notations. From {\bf C.6} and the definitions, it is immediate that the normed representation in the inertial class of $1$ is $(-1)^{n-1}$, if $H$ is unramified.

\null {\bf Proposition:} \it Assume that $H$ is unramified.  Denote by $\widehat{T}$ the Langlands dual of the maximal torus of $H$. Let $s$ be in $\widehat{T}$ such that $(-1)^{n-1}_{s,E}$ is a conjugate-dual representation of type $^LH$.

Write $s=diag(x_1,\dots ,x_{[{n\over 2}]},\widehat{1}, \ol{x}_{[{n\over 2}]}^{-1},\dots , \ol{x}_1^{-1})\in GL_n(\Bbb C)$ (with $1$ appearing only when $n$ is odd and $[{n\over 2}]$ denoting the integer part of ${n\over 2}$). For $x\in\{x_1,\dots ,x_{[{n\over 2}]}\}$, denote by $m(x,s)$ the multiplicity of $x$ in $s$ and put
$$C_x=\cases GL_{m(x;s),} & \hbox{\rm if $x\not\in\{\pm 1\}$}\cr O_{m(1,s),}  & \hbox{\rm if $x=1$,} \cr  Sp_{m(-1,s),}& \hbox{\rm if $x=-1$.} \cr\endcases$$

Then, $C_{\widehat{H}}(Im((-1)^n_s))$ is isomorphic to $\prod_x C_x$, the product going over equivalence classes of elements in the set $\{x_1,\dots ,x_{[{n\over 2}]}\}$ with respect to the relation $x\sim x^{-1}$.

\null Proof: \rm This follows from \cite{GGP, 8.1(iii) and section 4}. \hfill{\fin 2}

\null{\bf C.8} Proposition {\bf 2.1} remains valid, after replacing the Langlands parameter $\varphi $ by $\varphi _E$ and self-dual by conjugate-dual in appropriate places {\cite{GGP, sections 4 and 8}. In the same spirit, one gets the generalization of {\bf 2.2 - 2.6} with ${C'}^+:=C'$.

\null{\bf C.9} With all these changes, the corollary {\bf 3.5} is valid for the quasi-split unitary group.

\null
\null

\enddocument